\documentclass{article}
\usepackage{amsfonts}
\usepackage{amsmath}
\usepackage{amsmath}
\usepackage{amssymb}

\newtheorem{theorem}{Theorem}

\newtheorem{corollary}[theorem]{Corollary}

\newtheorem{lemma}[theorem]{Lemma}

\begin{document}

\title{Decomposable branching processes having a fixed extinction moment%
\thanks{%
This work is supported by the RSF under a grant 14-50-00005. } }
\author{Vatutin V.A.\thanks{%
Department of Discrete Mathematics, Steklov Mathematical Institute, 8,
Gubkin str., 119991, Moscow, Russia; e-mail: vatutin@mi.ras.ru}, Dyakonova
E.E.\thanks{%
Department of Discrete Mathematics, Steklov Mathematical Institute, 8,
Gubkin str., 119991, Moscow, Russia; e-mail: elena@mi.ras.ru}}
\maketitle

\begin{abstract}
The asymptotic behavior, as $n\rightarrow \infty $ of the probability of the
event that a decomposable critical branching process $\mathbf{Z}%
(m)=(Z_{1}(m),...,Z_{N}(m)),$ $m=0,1,2,...,$ with $N$ types of particles
dies at moment $n$ is investigated and conditional limit theorems are proved describing the
distribution of the number of particles in the process $\mathbf{Z}(\cdot)$
at moment $m<n,$ given that the extinction moment of the process is $n$.

These limit theorems may be considered as the statements describing the
distribution of the number of vertices in the layers of certain classes of
simply generated random trees having a fixed hight.
\end{abstract}



\section{Introduction}

We consider a Galton-Watson branching process with $N$ types of particles
labelled $1,2,...,N$ \ and denote by
\begin{equation*}
\mathbf{Z}(n)=\left( Z_{1}(n),...,Z_{N}\left( n\right) \right)
\end{equation*}%
the population vector at time $n\in \mathbb{Z}_{+}=\left\{ 0,1,...\right\} ,%
\mathbf{Z}(0)=\left( 1,0,...,0\right) .$ Denote by $T_{N}$ the extinction
moment of the process. The aim of the present paper is to investigate the
asymptotic behavior, as $n\rightarrow \infty $ of the probability of the
event $\left\{ T_{N}=n\right\} $ and the distribution of the random vector $%
\mathbf{Z}(m),0\leq m<n,$ given $T_{N}=n$ and assuming that $\mathbf{Z}%
(\cdot )$ is a decomposable critical branching process.

Properties of the single-type critical Galton-Watson process given its
extinction moment have been investigated by a number of authors (see, for
instance, \cite{KN93}, \cite{Kes86}, \cite{Pak71}). Asymptotic properties of
the survival probability for multitype indecomposable critical Markov
processes as well as the properties of these processes given their survival
up to a distant moment were analysed in \cite{JS67}, \cite{Mu63} and \cite%
{VV77}.

The decomposable branching processes are less investigated. We mention
papers \cite{FN}, \cite{FN2}, \cite{Og}, \cite{P76}, \cite{P77}, \cite%
{ChisSav}, \cite{VV14}, \cite{VS}, \cite{VDJS},\ \cite{Z} in this connection
where the asymptotic representations for the
probability of the event $\left\{ T_{N}>n\right\} $ are found under various restrictions and the
Yaglom-type limit theorems for the distribution of the number of particles
are proved for the multi-type decomposable critical Markov processes (and
their reduced analogues) under the condition $T_{N}>n.$ However, the study
of the asymptotic properties of the probability $\mathbf{P}\left(
T_{N}=n\right) $ for the decomposable critical Markov branching processes and
investigation of the conditional distributions of the number of particles in
these processes given $T_{N}=n,$ have not been considered up to now. The present paper
deals with such circle of questions.

Namely, we consider a decomposable Galton-Watson branching process with $N$
types of particles in which a type~$i$ parent-particle may produce children
of types $j\geq i$ only.

Introduce the probability generating functions for the distribution laws of
the offspring sizes of particles
\begin{equation}
h^{(i,N)}(\mathbf{s})=h^{(i,N)}(s_{i},...,s_{N})=\mathbf{E}\left[
s_{i}^{\eta _{i,i}}...\,s_{N}^{\eta _{i,N}}\right] ,\ i=1,2,...,N,
\label{DefNONimmigr}
\end{equation}
where the random variable $\eta _{i,j}$ is equal to the number of type~$j$
 daughter particles of a type~$i$ particle.

Let $\mathbf{e}_{i}$ be an $N$-dimensional vector whose $i$-th component is
equal to one while the remaining are zeros and $\mathbf{0}=(0,...,0)$ be an $%
N$--dimensional vector all whose components are zeros. The first moments of
the components of $\mathbf{Z}\left( n\right) $ will be denoted as%
\begin{equation*}
m_{i,j}(n)=\mathbf{E}\left[ Z_{j}\left( n\right) |\mathbf{Z}_{0}=\mathbf{e}%
_{i}\right]
\end{equation*}%
with $m_{i,j}=m_{i,j}(1)=\mathbf{E}[\eta_{i,j}]$ \thinspace being the
average number of type $j$ children produced by a particle of type $i$. \

Since $m_{i,j}=0$ if $i>j$, the mean matrix $\mathbf{M}%
=(m_{i,j})_{i,j=1}^{N} $ of our decomposable Galton-Watson branching process
has the form
\begin{equation}
\mathbf{M=}\left(
\begin{array}{ccccc}
m_{1,1} & m_{1,2} & ... & ... & m_{1,N} \\
0 & m_{2,2} & ... & ... & m_{2,N} \\
0 & 0 & m_{3,3} & ... & ... \\
... & ... & ... & ... & ... \\
... & ... & ... & ... & ... \\
0 & 0 & ... & 0 & m_{N,N}%
\end{array}%
\right) .  \label{matrix1}
\end{equation}

We say that \textbf{Hypothesis A} is valid if the decomposable branching
process with $N$ types of particles is strongly critical, i.e. (see \cite%
{FN2})
\begin{equation}
m_{i,i}=\mathbf{E}\left[ \eta _{i,i}\right] =1,\quad i=1,2,...,N
\label{Matpos}
\end{equation}%
and, in addition,
\begin{equation}
m_{i,i+1}=\mathbf{E}\left[ \eta _{i,i+1}\right] \in (0,\infty ),\
i=1,2,...,N-1,  \label{Maseq}
\end{equation}%
\begin{equation}
\mathbf{E}\left[ \eta _{i,j}\eta _{i,k}\right] <\infty ,\,i=1,...,N;\
k,j=i,i+1,...,N,  \label{FinCovar}
\end{equation}%
with
\begin{equation}
b_{i}=\frac{1}{2}Var\left[ \eta _{i,i}\right] \in \left( 0,\infty \right)
,i=1,2,...,N.  \label{FinVar}
\end{equation}%
Thus, a particle of the process is able to produce the direct descendants of
its own type, of the next in the order type, and (not necessarily, as direct
descendants) of all the remaining in the order types, but not any preceding
ones.

In the sequel we assume (if otherwise is not stated) that $\mathbf{Z}(0)=%
\mathbf{e}_{1}$, i.e. we suppose that the branching process in question is
initiated at time $n=0$ by a single particle of type 1.

The functions $\Phi _{i}=\Phi _{i}(\lambda _{i},\lambda _{i+1},...,\lambda
_{N}),\ i=1,2,...,N-1,$ being for $\lambda _{j}\geq 0,j=1,2,...,N$ solutions
of the equations
\begin{equation}
\sum_{k=i}^{N}\left( k-i+1\right) \lambda _{k}\frac{\partial \Phi _{i}}{%
\partial \lambda _{k}}=-b_{i}\Phi _{i}^{2}+\Phi
_{i}+\sum_{k=i}^{N}f_{k,i}\lambda _{k},\ i=1,2,...,N-1,  \label{GlobalDif}
\end{equation}%
with the initial conditions
\begin{equation*}
\Phi _{i}(\mathbf{0})=0,\ \frac{\partial \Phi _{i}(\mathbf{0})}{\partial
\lambda _{i}}=1
\end{equation*}%
and, for $k>i$
\begin{equation*}
\frac{\partial \Phi _{i}(\mathbf{0})}{\partial \lambda _{k}}=\frac{f_{k,i}}{%
k-i}=\frac{1}{(k-i)!}\prod_{j=i}^{k-1}m_{j,j+1}
\end{equation*}%
are important in the statements of the theorems to follow. Existence and
uniqueness of the solutions of the mentioned equations are established in
\cite{FN2}. Note that if $N=2,$ then
\begin{equation}
\Phi _{2}(\lambda _{1},\lambda _{2})=\sqrt{\frac{m_{1,2}\lambda _{2}}{b_{1}}}%
\frac{b_{1}\lambda _{1}+\sqrt{b_{1}m_{1,2}\lambda _{2}}\tanh \sqrt{%
b_{1}m_{1,2}\lambda _{2}}}{b_{1}\lambda _{1}\tanh \sqrt{b_{1}m_{1,2}\lambda
_{2}}+\sqrt{b_{1}m_{1,2}\lambda _{2}}}.  \label{DefSimpl}
\end{equation}

Let
\begin{equation}
c_{i,N}=\left( \frac{1}{b_{N}}\right) ^{1/2^{N-i}}\prod_{j=i}^{N-1}\left(
\frac{m_{j,j+1}}{b_{j}}\right) ^{1/2^{j-i+1}},  \label{Const2}
\end{equation}

\begin{equation}
D_{i}=(b_{i}m_{i,i+1})^{1/2^{i}}c_{1,i},i=1,2,...,N.  \label{DcConnection}
\end{equation}
Denote
\begin{equation*}
T_{ki}=\min \left\{ n\geq 1:Z_{k}(n)+Z_{k+1}(n)+...+Z_{i}(n)=0|\mathbf{Z}(0)=%
\mathbf{e}_{k}\right\}
\end{equation*}
the extinction moment of the population consisting of the particles of types
$k,k+1,...,i,$ given that the process was initiated at moment $n=0$ by a
single particle of type~$k$. To simplify notation, set $T_{i}=T_{1i}.$

We fix $N\geq 2$ and use, when it is necessary, the notation
\begin{equation*}
\gamma _{0}=0,\ \gamma _{i}=\gamma _{i}(N)=2^{-(N-i)},\ i=1,2,...,N.
\end{equation*}
Besides, we write $a(n)\sim b(n)$ if $\lim_{n\to\infty}a(n)/b(n)=1$ and $a(n)\ll b(n)$ if $\lim_{n\to\infty}a(n)/b(n)=0.$

Asymptotic properties of the probability that a critical decomposable
Galton-Watson branching process dies out at a fixed moment are described by
the following theorem.

\begin{theorem}
\label{T_loc}If Hypothesis $A$ is valid, then
\begin{equation*}
\mathbf{P}\left( T_{iN}=n\right) \sim \frac{g_{i,N}}{n^{1+\gamma _{i}}}%
,\quad i=1,2,...,N,
\end{equation*}
where
\begin{equation*}
g_{i,N}=\gamma _{i}c_{i,N}.
\end{equation*}
\end{theorem}

We now formulate four more theorems in which, given $T_{N}=n$ the limiting
distributions of the number of particles at moment $m$ are found depending
on the ratio between the parameters $m$ and $n.$

\begin{theorem}
\label{T_initial}If $n^{\gamma _{1}}\gg m\rightarrow \infty $, then
\begin{equation*}
\lim_{m\rightarrow \infty }\mathbf{E}\left[ \exp \left\{
-\sum_{l=1}^{N}\lambda _{l}\frac{Z_{l}(m)}{m^{l}}\right\} \Big|\,T_{N}=n%
\right] =\frac{\partial \Phi _{1}(\lambda _{1},\lambda _{2},...,\lambda _{N})%
}{\partial \lambda _{1}}.
\end{equation*}
\end{theorem}

We see that, given $\{T_{N}=n\}$ particles of all types present in the
process at the initial stage of its evolution.

\begin{theorem}
\label{T_sharp2}If $m\sim yn^{\gamma _{i}}$ for some $y>0$ and $i\in \left\{
1,2,...,N-1\right\} ,$ then, for any $s_{k}\in \left[ 0,1\right]
,k=1,...,i-1,$ and $\lambda _{l}\geq 0,l=i,...,N$
\begin{eqnarray*}
&&\lim_{m\rightarrow \infty }\mathbf{E}\left[ s_{1}^{Z_{1}(m)}\cdot \cdot
\cdot s_{i-1}^{Z_{i-1}(m)}\exp \left\{ -\sum_{l=i}^{N}\lambda _{l}\frac{%
Z_{l}(m)}{n^{\left( l-i+1\right) \gamma _{i}}}\right\} \Big|\,T_{N}=n\right]
\\
&&\quad =D_{i-1}\frac{g_{i,N}}{g_{1,N}}\frac{\partial }{\partial \lambda _{i}%
}\left( \frac{\Phi _{i}(\lambda _{i}^{\prime }y,\lambda _{i+1}^{\prime
}y^{2},\lambda _{i+2}y^{3},...,\lambda _{N}y^{N-i+1})}{y}\right) ^{1/2^{i}}
\\
&&\qquad +D_{i-1}\frac{g_{i+1,N}}{g_{1,N}}\frac{\partial }{\partial \lambda
_{i+1}}\left( \frac{\Phi _{i}(\lambda _{i}^{\prime }y,\lambda _{i+1}^{\prime
}y^{2},\lambda _{i+2}y^{3},...,\lambda _{N}y^{N-i+1})}{y}\right) ^{1/2^{i}},
\end{eqnarray*}
where $\lambda _{i}^{\prime }=\lambda _{i}+c_{i,N},\lambda _{i+1}^{\prime
}=\lambda _{i+1}+c_{i+1,N}.$
\end{theorem}

Observe that if $i>1,$ then, under the conditions of Theorem \ref{T_sharp2}
there are no particles of the types $1,2,...,i-1$ in the limit.

Let $a_{i,i}=1$ and for $i<j$
\begin{equation}
a_{i,j}=\frac{1}{\left( j-i\right) !}\prod_{k=i}^{j-1}m_{k,k+1}.  \label{len}
\end{equation}

\begin{theorem}
\label{T_interm2}If $n^{\gamma _{i}}\ll m\ll n^{\gamma _{i+1}}$ for some $%
i\in \left\{ 1,2,...,N-1\right\} ,$ then, for any $s_{k},k=1,2,...,i$ and $%
\lambda _{l}\geq 0,l=i+1,...,N$
\begin{eqnarray*}
&&\lim_{m\rightarrow \infty }\mathbf{E}\left[ s_{1}^{Z_{1}(m)}\cdot \cdot
\cdot s_{i}^{Z_{i}(m)}\exp \left\{ -\sum_{l=i+1}^{N}\lambda _{l}\frac{%
Z_{l}(m)}{n^{\gamma _{i+1}}m^{l-i-1}}\right\} \Big|\,T_{N}=n\right] \\
&&\qquad\qquad\qquad\qquad\qquad=\frac{D_{i}}{2^{i}}\frac{g_{i+1,N}}{g_{1,N}}%
\left( c_{i+1,N}+\sum_{l=i+1}^{N}\lambda _{l}a_{i+1,l}\right) ^{-1+1/2^{i}}.
\end{eqnarray*}
\end{theorem}

We see that, under the conditions of Theorem \ref{T_interm2} there are no
particles of types $1,2,...,i$ in the limit.

\begin{theorem}
\label{T_finalstage}If $m\sim xn,x\in \left( 0,1\right) ,$ then
\begin{eqnarray*}
&&\lim_{m\rightarrow \infty }\mathbf{E}\left[ s_{1}^{Z_{1}(m)}\cdot \cdot
\cdot s_{N-1}^{Z_{N-1}(m)}\exp \left\{ -\lambda _{N}\frac{Z_{N}(m)}{b_{N}n}%
\right\} \Big|\,T_{N}=n\right]  \\
&&\qquad \qquad \qquad =\frac{1}{\left( 1+(1-x)\lambda _{N}\right)
^{1-\gamma_1}}\frac{1}{\left( 1+\lambda _{N}x\left( 1-x\right) \right)
^{1+\gamma_1}}.
\end{eqnarray*}
\end{theorem}

It follows from Theorem \ref{T_finalstage} that at the final stage of the
development the population contains particles of type $N$ only.

We note that Theorems \ref{T_initial}-\ref{T_finalstage} may be considered
as the statements describing the distribution of the number of vertices in
the layers of certain classes simply generated random trees having a fixed
hight (see \cite{MM77}). The vertices of such trees are colored by one of $N$
colors labelled by numbers $1$ through $N,$ and the numbers of the colors
are monotone decreasing from the leaves to the root. The reader may find a
more detailed information about the properties of simply generated trees and
their connection with branching processes in a recent survey~\cite{Jan12}.

\section{Preliminary arguments}

In the sequel we denote by $\varepsilon _{i}(n), \varepsilon _{i}(n; m),
i=1,2,...$ some functions vanishing as $n\rightarrow \infty.$ These function
may be not necessary the same in different formulas.

\begin{lemma}
\label{L_basic}Let $A,B,\alpha$ and $\beta $ be positive numbers, $\alpha >\beta ,\beta
\in (0,1)$, and let $\Delta _{n},n=0,1,2,...$ be a sequence nonnegative
numbers meeting the recurrent relationships
\begin{equation}
\Delta _{0}=0,\Delta _{n}=\frac{A}{n^{\alpha }}\left( 1+\varepsilon
_{1}(n)\right) +\Delta _{n-1}\left( 1-\frac{B}{n^{\beta }}\left(
1+\varepsilon _{2}(n)\right) \right) ,\quad n=1,2,....  \label{recur}
\end{equation}
Then
\begin{equation*}
\lim_{n\rightarrow \infty }n^{\alpha -\beta }\Delta _{n}=\frac{A}{B}.
\end{equation*}
\end{lemma}

\textbf{Proof}. Set $\gamma =\alpha -\beta $ and write $\Delta _{n},n\geq 2$
in the form
\begin{equation}
\Delta _{n}=\frac{A}{B}\frac{1}{n^{\gamma }}\left( 1+\frac{\psi \left(
n\right) }{\log n}\right) .  \label{Reper}
\end{equation}%
Since
\begin{equation*}
\frac{1}{\left( n-1\right) ^{\gamma }}=\frac{1}{n^{\gamma }}\left( 1+\frac{%
\gamma }{n}\left( 1+\varepsilon _{3}(n)\right) \right) ,
\end{equation*}%
and
\begin{equation}
\frac{1}{\log \left( n-1\right) }=\frac{1}{\log n}\left( 1+\frac{1}{n\log n}%
\left( 1+\varepsilon _{4}(n)\right) \right) ,  \label{llog}
\end{equation}%
(\ref{recur}) takes the from
\begin{eqnarray*}
&&\frac{A}{B}\frac{1}{n^{\gamma }}\left( 1+\frac{\psi \left( n\right) }{\log
n}\right) =\frac{A}{n^{\alpha }}\left( 1+\varepsilon _{1}(n)\right)  \\
&&\qquad +\frac{A}{B}\frac{1}{n^{\gamma }}\left( 1+\frac{\psi \left(
n-1\right) }{\log (n-1)}\right) \left( 1+\frac{\gamma }{n}\left(
1+\varepsilon _{3}(n)\right) \right) \left( 1-\frac{B}{n^{\beta }}\left(
1+\varepsilon _{2}(n)\right) \right) ,
\end{eqnarray*}%
which, after evident transformations based on the condition $\beta <1$ and
the equalities
\begin{eqnarray*}
&&\frac{A}{B}\frac{1}{n^{\gamma }}\left( 1+\frac{\gamma }{n}\left(
1+\varepsilon _{3}(n)\right) \right) \left( 1-\frac{B}{n^{\beta }}\left(
1+\varepsilon _{2}(n)\right) \right)  \\
&&\qquad \qquad \qquad =\frac{A}{B}\frac{1}{n^{\gamma }}\left( 1-\frac{B}{%
n^{\beta }}\left( 1+\varepsilon _{4}(n)\right) \right) =\frac{A}{B}\frac{1}{%
n^{\gamma }}-\frac{A}{n^{\alpha }}+\frac{\varepsilon _{5}(n)}{n^{\alpha }},
\end{eqnarray*}%
leads to
\begin{equation*}
\frac{\psi \left( n\right) }{\log n}=\frac{\varepsilon _{5}(n)}{n^{\beta }}+%
\frac{\psi \left( n-1\right) }{\log (n-1)}\left( 1-\frac{B}{n^{\beta }}%
\left( 1+\varepsilon _{6}(n)\right) \right)
\end{equation*}%
or, on account of (\ref{llog}),
\begin{equation*}
\psi \left( n\right) =\frac{\varepsilon _{5}(n)\log n}{n^{\beta }}+\psi
\left( n-1\right) \left( 1-\frac{B}{n^{\beta }}\left( 1+\varepsilon
_{6}(n)\right) \right) .
\end{equation*}%
We show that
\begin{equation}
\lim \sup_{n\rightarrow \infty }\frac{\left\vert \psi \left( n\right)
\right\vert }{\log n}=0.  \label{llog2}
\end{equation}%
If this is not the case, then there exists such a subsequence $%
n_{k}\rightarrow \infty $ as $k\rightarrow \infty ,$ that
\begin{equation}
\lim \sup_{k\rightarrow \infty }\frac{\left\vert \psi \left( n_{k}\right)
\right\vert }{\log n_{k}}=c>0.  \label{Contr}
\end{equation}%
Assume that $\psi \left( n_{k}\right) \rightarrow \infty $. Let $k$ be such
that
\begin{equation*}
\psi \left( n_{k}\right) =\max_{1\leq n\leq n_{k}}\psi \left( n\right) \text{%
.}
\end{equation*}%
Then (to simplify notation we agree to write $n_{k}=n)$
\begin{equation*}
\psi \left( n\right) \leq \frac{\varepsilon _{5}(n)\log n}{n^{\beta }}+\psi
\left( n\right) \left( 1-\frac{B}{n^{\beta }}\left( 1+\varepsilon
_{6}(n)\right) \right)
\end{equation*}%
or
\begin{equation}
B\left( 1+\varepsilon _{6}(n)\right) \psi \left( n\right) \leq \varepsilon
_{5}(n)\log n.  \label{Above}
\end{equation}

Assume now that $\psi \left( n_{k}\right) \rightarrow -\infty $. Let $k$ be
such that
\begin{equation*}
\psi \left( n_{k}\right) =\min_{1\leq n\leq n_{k}}\psi \left( n\right) \text{%
.}
\end{equation*}%
Then (to simplify notation we agree to write $n_{k}=n)$
\begin{equation*}
\psi \left( n\right) \geq \frac{\varepsilon _{5}(n)\log n}{n^{\beta }}+\psi
\left( n\right) \left( 1-\frac{B}{n^{\beta }}\left( 1+\varepsilon
_{6}(n)\right) \right)
\end{equation*}%
or
\begin{equation}
B\left( 1+\varepsilon _{6}(n)\right) \psi \left( n\right) \geq \varepsilon
_{5}(n)\log n.  \label{below}
\end{equation}%
Clearly, the combination of (\ref{Above}) and (\ref{below}) contradicts (\ref%
{Contr}). This proves (\ref{llog2}).

It follows from the obtained estimate and (\ref{Reper}) that
\begin{equation*}
\Delta _{n}\sim \frac{A}{B}\frac{1}{n^{\alpha -\beta }},\quad n\rightarrow
\infty .
\end{equation*}%
Lemma \ref{L_basic} is proved.

We use the symbols $\mathbf{P}_{i}$ and $\mathbf{E}_{i}$ to denote the
probability and expectation calculated under the condition that a branching
process is initiated at moment $n=0$ by a single particle of type $i$.
Sometimes we write $\mathbf{P}$ and $\mathbf{E}$ for $\mathbf{P}_{1}$ and $%
\mathbf{E}_{1},$ respectively.

Denote by
\begin{equation}
b_{ikl}(n)=\mathbf{E}_{i}\left[ Z_{k}(n)Z_{l}(n)-\delta _{kl}Z_{l}(n)\right]
\label{Defq}
\end{equation}%
the second moments of the components of the process $\mathbf{Z}\left(
n\right) $. Let$\ b_{ikl}=b_{ikl}(1).$

For any vector $\mathbf{s}=(s_{1},...,s_{p})$ (the dimension will usually be
clear from the context), and any vector $\mathbf{k}=(k_{1}.....k_{p})$ $\ $%
with integer valued components define%
\begin{equation*}
\mathbf{s}^{\mathbf{k}}=s_{1}^{k_{1}}...\,s_{p}^{k_{p}}.
\end{equation*}%
Further, let $\mathbf{1}=\left( 1,...,1\right) $ be a vector of units.
Sometimes it will be convenient to write $\mathbf{1}^{(i)}$ for the $i-$%
dimensional vector with all components equal to one.

Let
\begin{equation*}
H_{n}^{(i,N)}(\mathbf{s})=\mathbf{E}_{i}\left[ \mathbf{s}^{\mathbf{Z}(n)}%
\right] =\mathbf{E}_i\left[ s_{i}^{Z_{i}(n)}...\,s_{N}^{Z_{N}(n)}\right]
\end{equation*}
be the probability generating functions for the process $\mathbf{Z}(n)$
given the process is initiated by a single particle of type $i\in \left\{
1,2,...,N\right\} $ at moment $0$. Clearly (see (\ref{DefNONimmigr})), $%
H_{1}^{(i,N)}(\mathbf{s})=h^{(i,N)}(\mathbf{s})$ for any $i\in \{1,...,N\}.$
Denote
\begin{equation*}
Q_{n}^{(i,N)}(\mathbf{s})=1-H_{n}^{(i,N)}(\mathbf{s}),\quad
Q_{n}^{(i,N)}=1-H_{n}^{(i,N)}(\mathbf{0})
\end{equation*}
and let
\begin{equation*}
\mathbf{H}_{n}(\mathbf{s})=(H_{n}^{(1,N)}(\mathbf{s}),...,H_{n}^{(N,N)}(%
\mathbf{s})),~\mathbf{Q}_{n}(\mathbf{s})=(Q_{n}^{(1,N)}(\mathbf{s}%
),...,Q_{n}^{(N,N)}(\mathbf{s})).
\end{equation*}

The starting point of our arguments is the following theorem being a
simplified combination of the respective results from \cite{FN} and \cite%
{FN2}:

\begin{theorem}
\label{T_Foster}Let $\mathbf{Z}(n),n=0,1,...$ be a decomposable critical
branching process meeting the conditions (\ref{matrix1}), (\ref{Matpos}), (%
\ref{Maseq}) and (\ref{FinCovar}). Then $m_{j,j}(n)=1$ and, as $n\rightarrow
\infty $
\begin{eqnarray}
m_{i,j}(n) &\sim &a_{i,j}n^{j-i},\ i<j,  \label{MomentSingle3} \\
b_{jpq}(n) &=&a_{jpq}n^{p+q-2j+1}+o\left( n^{p+q-2j+1}\right) ,\ j\leq \min
(p,q),  \label{Momvariance}
\end{eqnarray}
where $a_{i,j}$ are the same as in (\ref{len}) and $a_{jpq}$ are nonnegative
constants known explicitly (see \cite{FN2}, Theorem 1).

In addition \textit{(}see \cite{FN}, Theorem~1\textit{)}, as $n\rightarrow
\infty $
\begin{equation}
Q_{n}^{(i,N)}=1-H_{n}^{(i,N)}(\mathbf{0})=\mathbf{P}_i(\mathbf{Z}(n)\neq
\mathbf{0})\sim c_{i,N}n^{-1/2^{N-i}},  \label{SurvivSingle}
\end{equation}
where $c_{i,N}$ are the same as in (\ref{Const2}), and for $\lambda _{l}\geq
0,\,l=1,2,...,N,$
\begin{eqnarray}
&&\lim_{m\rightarrow \infty }m\left( 1-\mathbf{E}\left[ \exp \left\{
-\sum_{l=1}^{N}\lambda _{l}\frac{Z_{l}(m)}{m^{l}}\right\} \right] \right)
\notag \\
&&\qquad \qquad =\lim_{m\rightarrow \infty }mQ_{m}^{(1,N)}\Big(e^{-\lambda
_{1}/m},e^{-\lambda _{2}/m^{2}},...,e^{-\lambda _{N}/m^{N}}\Big)  \notag \\
&&\qquad \qquad =\Phi _{1}(\lambda _{1},\lambda _{2},...,\lambda _{N}).
\label{Fost1}
\end{eqnarray}
\end{theorem}

We need also the following Yaglom-type limit theorem proved in \cite{VV14}
and complementing Theorem \ref{T_Foster}.

\begin{theorem}
\label{T_Yaglom}If the conditions of Theorem \ref{T_Foster} are valid, then
for any $\lambda >0$
\begin{equation}
\lim_{n\rightarrow \infty }\mathbf{E}_{1}\left[ \exp \left\{ -\lambda \frac{%
Z_{N}(n)}{b_{N}n}\right\} \Big|\mathbf{Z}(n)\neq \mathbf{0}\right] =1-\Big(%
\frac{\lambda }{1+\lambda }\Big)^{1/2^{N-1}}.  \label{Yag}
\end{equation}
\end{theorem}

\section{Proof of Theorem \protect\ref{T_loc}}

The results of the previous section allow us to prove Theorem \ref{T_loc}.
For $N=i$ we have
\begin{eqnarray*}
\mathbf{P}\left( T_{NN}=n\right) &=&H_{n}^{(N,N)}(0)-H_{n-1}^{(N,N)}(0) \\
&=&h^{(N,N)}(H_{n-1}^{(N,N)}(0))-h^{(N,N)}(H_{n-2}^{(N,N)}(0)) \\
&\leq &H_{n-1}^{(N,N)}(0)-H_{n-2}^{(N,N)}(0)=\mathbf{P}\left(
T_{NN}=n-1\right) .
\end{eqnarray*}%
Since
\begin{equation}
\mathbf{P}\left( T_{NN}>n\right) =\sum_{k=n+1}^{\infty }\mathbf{P}\left(
T_{NN}=k\right) \sim \frac{1}{b_{N}n}  \label{Loc1}
\end{equation}%
as $n\rightarrow \infty $ and the sequence $\mathbf{P}\left( T_{NN}=k\right)
$ is monotone, (\ref{Loc1}) and Corollary 2 in~\cite{VV77a} imply as $%
n\rightarrow \infty $
\begin{equation*}
\mathbf{P}\left( T_{NN}=n\right) \sim \frac{1}{b_{N}n^{2}},
\end{equation*}%
proving Theorem \ref{T_loc} for $i=N$. Assume that the theorem is proved for all $%
i\in \left\{ j+1,N\right\} $, where $1<j+1\leq N$. Let us demonstrate that
it is true for $i=j$.

To this aim we put

\begin{equation*}
\varphi (t)=h^{(j,N)}\left( \mathbf{H}_{n-2}(\mathbf{0})+t(\mathbf{H}_{n-1}(%
\mathbf{0})-\mathbf{H}_{n-2}(\mathbf{0}))\right),\, 0\leq t\leq 1.
\end{equation*}
Clearly that

\begin{equation}
\mathbf{P}\left( T_{jN}=n\right) =\varphi (1)-\varphi (0)=\varphi ^{\prime
}(0)+\varphi ^{\prime \prime }(\theta _{n})/2,  \label{PP}
\end{equation}%
where $\theta _{n}\in \lbrack 0,1].$

It is easy to check that
\begin{eqnarray*}
\varphi^{\prime} (0)&=&\sum_{i=j}^N\frac{\partial h^{(j,N)}(\mathbf{H}_{n-2}(%
\mathbf{0}))}{\partial s_i}(H^{(i,N)}_{n-1}(\mathbf{0})-H^{(i,N)}_{n-2}(%
\mathbf{0})) \\
&=&(1+\varepsilon (n))\sum_{i=j+1}^Nm_{j,i}\mathbf{P}\left(
T_{iN}=n-1\right)+\frac{\partial h^{(j,N)}(\mathbf{H}_{n-2}(\mathbf{0}))}{%
\partial s_j}\mathbf{P}\left( T_{jN}=n-1\right),
\end{eqnarray*}
where by the induction assumption
\begin{equation*}
\sum_{i=j+1}^Nm_{j,i}\mathbf{P}\left( T_{iN}=n-1\right)=(1+\varepsilon_1(n))%
\frac{m_{j,j+1}g_{j+1,N}}{n^{1+\gamma_{j+1}}},
\end{equation*}
and, in view of (\ref{SurvivSingle})
\begin{eqnarray*}
\frac{\partial h^{(j,N)}(\mathbf{H}_{n-2}(\mathbf{0}))}{\partial s_j}%
&=&1-(1+\varepsilon_2(n))\sum_{k=j}^N\frac{\partial^{2} h^{(j,N)}(\mathbf{1})%
}{\partial s_k\partial s_j}(1-H^{(k,N)}_{n-2}(\mathbf{0})) \\
&=&1-(1+\varepsilon_3(n))b_{jjj}(1-H^{(j,N)}_{n-2}(\mathbf{0}%
))=1-(1+\varepsilon_4(n))\frac{c_{j,N}b_{jjj}}{n^{\gamma_j}}.
\end{eqnarray*}

Further, for $\Theta _{n}=\mathbf{H}_{n-2}(\mathbf{0})+\theta _{n}(\mathbf{H}%
_{n-1}(\mathbf{0})-\mathbf{H}_{n-2}(\mathbf{0}))$ we have
\begin{eqnarray*}
\varphi ^{\prime \prime }(\theta _{n}) &=&\sum_{k=j}^{N}\sum_{i=j}^{N}\frac{%
\partial h^{(j,N)}(\Theta _{n})}{\partial s_{k}\partial s_{i}}%
(H_{n-1}^{(k,N)}(\mathbf{0})-H_{n-2}^{(k,N)}(\mathbf{0}))(H_{n-1}^{(i,N)}(%
\mathbf{0})-H_{n-2}^{(i,N)}(\mathbf{0})) \\
&=&(1+\varepsilon _{5}(n))\sum_{k=j}^{N}\sum_{i=j}^{N}b_{jik}\mathbf{P}%
\left( T_{kN}=n-1\right) \mathbf{P}\left( T_{iN}=n-1\right) \\
&=&(1+\varepsilon _{5}(n))b_{jjj}\mathbf{P}^{2}\left( T_{j,N}=n-1\right) +o%
\Big(\frac{1}{n^{1+\gamma _{j+1}}}\Big).
\end{eqnarray*}

Substituting the obtained estimates in (\ref{PP}) and recalling that $%
b_{j}=b_{jjj}/2,$ we get

\begin{eqnarray*}
\mathbf{P}\left( T_{jN}=n\right) &=&\frac{m_{j,j+1}g_{j+1,N}}{n^{1+\gamma
_{j+1}}}\left( 1+\varepsilon _{1}(n)\right) \\
&&+\mathbf{P}\left( T_{jN}=n-1\right) \left( 1-\frac{2b_{j}c_{j,N}}{%
n^{\gamma _{j}}}\left( 1+\varepsilon _{2}(n)\right) \right) .
\end{eqnarray*}%
This representation, Lemma \ref{L_basic} and the equalities
\begin{equation*}
\frac{m_{j,j+1}g_{j+1,N}}{2b_{j}}=\gamma _{j+1}\frac{m_{j,j+1}c_{j+1,N}}{%
2b_{j}}=\gamma _{j}c_{j,N}^{2}=c_{j,N}g_{jN}
\end{equation*}%
yield, as $n\rightarrow \infty $
\begin{equation*}
\mathbf{P}\left( T_{jN}=n\right) \sim \frac{m_{j,j+1}g_{j+1,N}}{2b_{j}c_{j,N}%
}\frac{1}{n^{1+\gamma _{j+1}-\gamma _{j}}}=\frac{g_{j,N}}{n^{1+\gamma _{j}}}.
\end{equation*}
This proves Theorem \ref{T_loc} by induction.

\section{Auxiliary lemmas}

We prove in this section a number of statements about the asymptotic
behavior, as $n\rightarrow \infty $ expectations of the form
\begin{equation*}
\mathbf{E}\left[ \exp \left\{ -\sum_{l=1}^{N}\lambda _{l}\frac{Z_{l}(m)}{%
r(n,m)}\right\} \right]
\end{equation*}
and their derivatives with respect to the parameters $\lambda
_{l},l=1,2,...,N$ depending on the rate of growth the parameter $m=m(n)$ to
infinity and the form of scaling $r(n,m)$. We will show that the asymptotic
behavior of the mentioned quantities is essentially different for the cases $%
m\ll n^{\gamma _{1}},m\sim yn^{\gamma _{i}},y>0,$ $n^{\gamma _{i}}\ll m\ll
n^{\gamma _{i+1}},i=1,2,...,N-1$ and $m\sim xn,x\in \left( 0,1\right) $.

\subsection{The case $m\sim yn^{\protect\gamma _{i}},y>0,1\leq i\leq N-1$}

Let%
\begin{equation*}
I_{k}(m)=I\{Z_{1}(m)+\cdots +Z_{k}(m)=0\}
\end{equation*}%
be the indicator of the event that  particles of types $1,2,...,k$ are\
absent in the population at time $m$. Suppose that $I_{0}(m)=1$.

The aim of the present subsection is to prove the following lemma.

\begin{lemma}
\label{L_sharp!} If the asymptotic relation $m\sim yn^{\gamma _{i}},y>0,$ is
true for some $i\in \{1,...,N-1\},$ then for any $%
j\in \left\{ i,...,N-1\right\} $ and any tuple $\lambda _{l}\geq 0,l=i,...,N$
\begin{eqnarray*}
&&\lim_{n\rightarrow \infty }\frac{n^{\gamma _{1}}}{n^{\left( j-i+1\right)
\gamma _{i}}}\mathbf{E}\left[ Z_{j}(m)\exp \left\{ -\sum_{l=i}^{N}\lambda
_{l}\frac{Z_{l}(m)}{n^{\left( l-i+1\right) \gamma _{i}}}\right\} I_{i-1}(m)%
\right] \\
&&\qquad =\lim_{n\rightarrow \infty }\frac{n^{\gamma _{1}}}{n^{\left(
j-i+1\right) \gamma _{i}}}\mathbf{E}\left[ Z_{j}(m)\exp \left\{
-\sum_{l=i}^{N}\lambda _{l}\frac{Z_{l}(m)}{n^{\left( l-i+1\right) \gamma
_{i}}}\right\} \right] \\
&&\qquad \qquad \qquad \quad =D_{i-1}\frac{\partial }{\partial \lambda _{j}}%
\left( \frac{\Phi _{i}(\lambda _{i}y,\lambda _{i+1}y^{2},...,\lambda
_{N}y^{N-i+1})}{y}\right) ^{1/2^{i-1}}.
\end{eqnarray*}
\end{lemma}

The desired statement will be a corollary of a number of lemmas the first of
them looks as follows.

\begin{lemma}
\label{L_twoOnly} If the asymptotic relation $m\sim yn^{\gamma _{i}},y>0,$
is true for some $i\in \{1,...,N-1\},$ then for $\lambda_l\geq 0,\, l=i,...,N$ and
\begin{equation*}
\mathbf{s}(i)=\left( \exp \left\{ -\frac{\lambda _{i}}{n^{\gamma _{i}}}%
\right\} ,\exp \left\{ -\frac{\lambda _{i+1}}{n^{2\gamma _{i}}}\right\}
,...,\exp \left\{ -\frac{\lambda _{N}}{n^{\left( N-i+1\right) \gamma _{i}}}%
\right\} \right)
\end{equation*}%
we have
\begin{equation*}
\lim_{n\rightarrow \infty }n^{\gamma _{i}}Q_{m}^{(i,N)}(\mathbf{s}%
(i))=y^{-1}\Phi _{i}(\lambda _{i}y,\lambda _{i+1}y^{2},...,\lambda
_{N}y^{N-i+1}).
\end{equation*}
\end{lemma}

\textbf{Proof.} Using (\ref{Fost1}) it is easy to check
\begin{eqnarray*}
&&\lim_{n\rightarrow \infty }n^{\gamma _{1}}Q_{m}^{(1,N)}(\mathbf{s}%
(i))=y^{-1}\lim_{m\rightarrow \infty }mQ_{m}^{(1,N)}(\mathbf{s}(i)) \\
&&\quad =y^{-1}\lim_{m\rightarrow \infty }m\left( 1-\mathbf{E}\left[ \exp
\left\{ -\sum_{l=1}^{N}\lambda _{l}\frac{Z_{l}(m)}{n^{l\gamma _{1}}}\right\} %
\right] \right)  \\
&&\quad =y^{-1}\lim_{m\rightarrow \infty }m\left( 1-\mathbf{E}\left[ \exp
\left\{ -\sum_{l=1}^{N}\lambda _{l}y^{l}\frac{Z_{l}(m)}{m^{l}}\right\} %
\right] \right)  \\
&&\quad =y^{-1}\Phi _{1}(\lambda _{1}y,\lambda _{2}y^{2},...,\lambda
_{N}y^{N}),
\end{eqnarray*}%
which proves the lemma for $i=1$. The cases when $i\in \{2,...,N-1\}$ may be
considered in a similar way.

Lemma \ref{L_twoOnly} is proved.

\begin{lemma}
\label{L_OnlyOneNew}If $m\sim yn^{\gamma _{j}},\,y>0,$  and $%
1\leq i<j\leq N,$ then for $\lambda_l\geq 0,\, l=j,...,N$
\begin{equation*}
\lim_{n\rightarrow \infty }n^{\left( j-i+1\right) \gamma _{i}}\left( 1-%
\mathbf{E}_{j}\left[ \exp \left\{ -\sum_{l=j}^{N}\lambda _{l}\frac{Z_{l}(m)}{%
n^{\left( l-i+1\right) \gamma _{i}}}\right\} \right] \right)
=\sum_{l=j}^{N}\lambda _{l}y^{l-j}a_{j,l}.
\end{equation*}
\end{lemma}

\textbf{Proof}. Set
\begin{equation*}
s_{l}=\exp \{-\frac{\lambda _{l}}{n^{(l-i+1)\gamma _{i}}}\},\,l=i,...,N,
\end{equation*}%
and consider the case $i=1$ only, since the proof for  $i\in \left\{
2,...,N-1\right\} $ requires only minor changes.

Clearly,
\begin{equation*}
0\leq \sum_{l=j}^{N}\left( 1-s_{l}\right) \mathbf{E}_{j}\left[ Z_{l}(m)%
\right]- Q_{m}^{(j,N)}(\mathbf{s})\leq \sum_{p,q=j}^{N}\left( 1-s_{p}\right)
\left( 1-s_{q}\right) \mathbf{E}_{j}\left[ Z_{p}(m)Z_{q}(m)\right] .
\end{equation*}

By (\ref{MomentSingle3}) we have as $n\rightarrow \infty $
\begin{eqnarray}
&&\sum_{l=j}^{N}(1-s_l)\mathbf{E}_{j}Z_{l}(m)\sim \sum_{l=j}^{N}\lambda _{l}%
\frac{\mathbf{E}_{j}Z_{l}(m)}{n^{l\gamma _{1}}}  \notag \\
&&\qquad\qquad\qquad\quad=\frac{1}{n^{j\gamma _{1}}}\sum_{l=j}^{N}\lambda
_{l}y^{l-j}\frac{\mathbf{E}_{j}Z_{l}(m)}{m^{l-j}}\sim \frac{1}{n^{j\gamma
_{1}}}\sum_{l=j}^{N}\lambda _{l}y^{l-j}a_{j,l}.  \label{add1}
\end{eqnarray}

Further, in view of (\ref{Momvariance})
\begin{eqnarray*}
&&\sum_{p,q=j}^{N}\lambda _{p}\lambda _{q}\frac{\mathbf{E}_{j}\left[
Z_{p}(m)Z_{q}(m)\right] }{n^{p\gamma _{1}}n^{q\gamma _{1}}} =\frac{1}{%
n^{j\gamma _{1}}}\sum_{p,q=j}^{N}\lambda _{p}\lambda _{q}\frac{b_{jpq}\left(
m\right) }{n^{\left( p-j\right) \gamma _{1}}n^{q\gamma _{1}}} \\
&&\quad=\frac{1}{n^{j\gamma _{1}}}\sum_{p,q=j}^{N}\lambda _{p}\lambda _{q}%
\frac{a_{jpq}m^{p+q-2j+1}}{n^{\left( p-j\right) \gamma _{1}}n^{q\gamma _{1}}}
+\varepsilon_1(n)\frac{1}{n^{j\gamma _{1}}}\sum_{p,q=j}^{N}\frac{m^{p+q-2j+1}%
}{n^{\left( p-j\right) \gamma _{1}}n^{q\gamma _{1}}} \\
&&\qquad\qquad\qquad\qquad\qquad\leq C\frac{1}{n^{j\gamma _{1}}}%
\sum_{p,q=j}^{N}\frac{n^{\gamma _{1}\left( p+q-2j+1\right) }}{n^{\left(
p-j\right) \gamma _{1}}n^{q\gamma _{1}}}=o\left( \frac{1}{n^{j\gamma _{1}}}%
\right) .
\end{eqnarray*}

The obtained estimates prove the lemma.

Let $\eta _{r,j}\left( k,l\right) $ be the number of type $j$ daughter
particles of the $l-$th particle of type~$r,$ belonging to the $k-$th
generation and let
\begin{equation*}
W_{p,i,j}=\sum_{r=p}^{i}\sum_{k=0}^{T_{i}}\sum_{q=1}^{Z_{r}(k)}\eta
_{r,j}\left( k,q\right)
\end{equation*}
be the total number of type $j\geq i+1$ daughter particles generated by all
the particles of types $p,p+1,...,i$ ever born in the process given that the
process is initiated at time $n=0$ by a single particle of type $p\leq i.$
Finally, put
\begin{equation*}
W_{p,i}=\sum_{j=i+1}^{N}W_{p,i,j}=\sum_{j=i+1}^{N}\sum_{r=p}^{i}%
\sum_{k=0}^{T_{i}}\sum_{q=1}^{Z_{r}(k)}\eta _{r,j}\left( k,q\right) .
\end{equation*}

\begin{lemma}
\label{L_Laplace} (see \cite{VV14}, Lemma 1). Let Hypothesis $A$ be valid.
Then, as $\lambda \downarrow 0$
\begin{equation}
1-\mathbf{E}\left[ e^{-\lambda W_{1,i,i+1}}\,|\mathbf{Z}(0)=\mathbf{e}_{1}%
\right] \sim D_{i}\lambda ^{1/2^{i}}  \label{Tot1}
\end{equation}
and there exists a constant $F_{i}>0$ such that
\begin{equation}
1-\mathbf{E}\left[ e^{-\lambda W_{1,i}}|\mathbf{Z}(0)=\mathbf{e}_{1}\right]
\sim F_{i}\lambda ^{1/2^{i}}.  \label{Tot2}
\end{equation}
\end{lemma}

Basing on Lemmas \ref{L_OnlyOneNew} and \ref{L_Laplace} we prove the
following statement.

\begin{lemma}
\label{L_initial} If $m\sim yn^{\gamma _{j}},\,y>0,$ for some $i\in \{1,2,...,N-1\},$ then for $\lambda_l\geq 0,\, l=i,...,N$
\begin{eqnarray}
&&\lim_{n\rightarrow \infty }n^{\gamma _{1}}\mathbf{E}\left[ \left( 1-\exp
\left\{ -\sum_{l=i}^{N}\lambda _{l}\frac{Z_{l}(m)}{n^{\left( l-i+1\right)
\gamma _{i}}}\right\} \right) I_{i-1}(m)\right]  \notag \\
&&\qquad =\lim_{n\rightarrow \infty }n^{\gamma _{1}}\mathbf{E}\left[ 1-\exp
\left\{ -\sum_{l=i}^{N}\lambda _{l}\frac{Z_{l}(m)}{n^{\left( l-i+1\right)
\gamma _{i}}}\right\} \right]  \notag \\
&&\qquad \qquad =D_{i-1}\left( \frac{\Phi _{i}(\lambda _{i}y,\lambda
_{i+1}y^{2},...,\lambda _{N}y^{N-i+1})}{y}\right) ^{1/2^{i-1}}.  \label{pre1}
\end{eqnarray}
\end{lemma}

\textbf{Proof.} For $i=1$ the statement of the lemma is a particular case of
Lemma~\ref{L_twoOnly}. Thus, we assume now that $i\in
\{2,3,...,N-1\}.$ According to~(\ref{SurvivSingle}) for $m\sim yn^{\gamma
_{i}}$ the following relations are valid:%
\begin{eqnarray*}
&&\lim_{n\rightarrow \infty }n^{\gamma _{1}}\mathbf{E}\left[ \left( 1-\exp
\left\{ -\sum_{l=i}^{N}\lambda _{l}\frac{Z_{l}(m)}{n^{\left( l-i+1\right)
\gamma _{i}}}\right\} \right) (1-I_{i-1}(m))\right] \\
&&\qquad \qquad \leq \lim_{n\rightarrow \infty }n^{\gamma _{1}}\mathbf{P}%
(T_{i-1}>m)=\lim_{n\rightarrow \infty }\frac{c_{1,i-1}n^{1/2^{N-1}}}{%
(yn^{1/2^{N-i}})^{1/2^{i-2}}}=0.
\end{eqnarray*}%
Therefore, to prove the lemma it is sufficient to show the validity of the
second equality in (\ref{pre1}) only. Recalling (\ref{SurvivSingle}) once more,
we have, as $n\rightarrow \infty ,$
\begin{equation*}
\mathbf{P}(T_{1,i-1}>n^{3\gamma _{i-2}})\sim c_{1,i-1}n^{-3\gamma
_{i-2}/2^{i-2}}=c_{1,i-1}n^{-3\gamma _{1}/2}=o(n^{-\gamma _{1}}).
\end{equation*}

Thus,
\begin{eqnarray*}
Q_{m}^{(1,N)}(\mathbf{s}) &=&\mathbf{E}\left[
1-s_{1}^{Z_{1}(m)}s_{2}^{Z_{2}(m)}...\,s_{N}^{Z_{N}(m)}\right] \\
&=&\mathbf{E}\left[ \left( 1-s_{i}^{Z_{i}(m)}...\,s_{N}^{Z_{N}(m)}\right)
;T_{i-1}\leq n^{3\gamma _{i-2}}\right] +o(n^{-\gamma _{1}}) \\
&=&1-H_{m}^{(1,N)}\left( \mathbf{1}^{(i-1)},s_{i},...,s_{N}\right)
+o(n^{-\gamma _{1}}).
\end{eqnarray*}%
It is not difficult to check that for our decomposable branching process
\begin{eqnarray*}
&&H_{m}^{(1,N)}\left( \mathbf{1}^{(i-1)},s_{i},...,s_{N}\right) \\
&&\quad =\mathbf{E}\left[ \prod_{k=0}^{m-1}\prod_{r=1}^{i-1}%
\prod_{l=1}^{Z_{r}(k)}\prod_{j=i}^{N}\left( H_{m-k}^{(j,N)}(\mathbf{s}%
)\right) ^{\eta _{r,j}\left( k,l\right) }\right] \\
&&\quad =\mathbf{E}\left[ \prod_{k=0}^{m-1}\prod_{r=1}^{i-1}%
\prod_{l=1}^{Z_{r}(k)}\prod_{j=i}^{N}\left( H_{m-k}^{(j,N)}(\mathbf{s}%
)\right) ^{\eta _{r,j}\left( k,l\right) };T_{i-1}\leq n^{3\gamma _{i-2}}%
\right] +o(n^{-\gamma _{1}}).\qquad
\end{eqnarray*}%
Observing that $\lim_{m\rightarrow \infty }H_{m-k}^{(j,N)}(\mathbf{s})=1$
for $k\leq T_{i-1}\leq n^{3\gamma _{i-2}}=o(m)$ and $j\geq i+1,$ we conclude
that, on the set $T_{i-1}\leq n^{3\gamma _{i-2}}$
\begin{eqnarray*}
&&\prod_{k=0}^{m-1}\prod_{r=1}^{i-1}\prod_{l=1}^{Z_{r}(k)}\prod_{j=i}^{N}%
\left( H_{m-k}^{(j,N)}(\mathbf{s})\right) ^{\eta _{r,j}\left( k,l\right) } \\
&&\quad =\exp \left\{
-\sum_{r=1}^{i-1}\sum_{k=0}^{T_{i-1}}\sum_{l=1}^{Z_{r}(k)}\sum_{j=i}^{N}\eta
_{r,j}\left( k,l\right) Q_{m-k}^{(j,N)}(\mathbf{s})(1+o(1))\right\} .
\end{eqnarray*}%
Note now that according to (\ref{Fost1}) for $m\sim yn^{\gamma _{i}},$
\begin{equation*}
s_{l}=\exp \left\{ -\frac{\lambda _{l}}{n^{\left( l-i+1\right) \gamma _{i}}}%
\right\} ,l=i,...,N
\end{equation*}%
and $k=o(m)$ the following relations are valid:
\begin{eqnarray*}
\lim_{n\rightarrow \infty }n^{\gamma _{i}}Q_{m-k}^{(i,N)}(\mathbf{s})
&=&y^{-1}\lim_{m\rightarrow \infty }m\left( 1-\mathbf{E}_{i}\left[ \exp
\left\{ -\sum_{l=i}^{N}\lambda _{l}y^{l-i+1}\frac{Z_{l}(m)}{m^{l-i+1}}%
\right\} \right] \right) \\
&=&y^{-1}\Phi _{i}(\lambda _{i}y,\lambda _{i+1}y^{2},...,\lambda
_{N}y^{N-i+1}),
\end{eqnarray*}%
while by Lemma \ref{L_OnlyOneNew} we have for $j>i$
\begin{equation}
\lim_{n\rightarrow \infty }n^{(j-i+1)\gamma _{i}}Q_{m-k}^{(j,N)}(\mathbf{s}%
)=\sum_{l=j}^{N}\lambda _{l}y^{l-j}a_{j,l}.  \label{oh}
\end{equation}%
Hence it follows that if the condition $T_{i-1}\leq \sqrt{mn^{\gamma _{i-1}}}%
=o(n^{\gamma _{i}})$ is valid, then
\begin{eqnarray*}
&&\sum_{r=1}^{i-1}\sum_{k=0}^{T_{i-1}}\sum_{l=1}^{Z_{r}(k)}\sum_{j=i}^{N}%
\eta _{r,j}\left( k,l\right) Q_{m-k}^{(j,N)}(\mathbf{s}) \\
&&\quad =(1+o(1))\sum_{j=i}^{N}Q_{m}^{(j,N)}(\mathbf{s})\sum_{r=1}^{i-1}%
\sum_{k=0}^{T_{i-1}}\sum_{l=1}^{Z_{r}(k)}\eta _{r,j}\left( k,l\right) \\
&&\quad =(1+o(1))\sum_{j=i}^{N}W_{1,i-1,j}Q_{m}^{(j,N)}(\mathbf{s}).
\end{eqnarray*}

Further, for $m\sim yn^{\gamma _{i}}$
\begin{equation}
W_{1,i-1,i}Q_{m}^{(i,N)}(\mathbf{s})=\left( 1+o(1)\right)
W_{1,i-1,i}y^{-1}\Phi _{i}(\lambda _{i}y,\lambda _{i+1}y^{2},...,\lambda
_{N}y^{N-i+1})n^{-\gamma _{i}},  \label{list}
\end{equation}%
while by (\ref{oh})
\begin{equation*}
\sum_{j=i+1}^{N}W_{1,i-1,j}Q_{m}^{(j,N)}(\mathbf{s})=O\left(
\sum_{j=i+1}^{N}W_{1,i-1,j}n^{-(j-i+1)\gamma _{i}}\right) =O(n^{-2\gamma
_{i}}W_{1,i-1}).
\end{equation*}%
Using Lemma \ref{L_Laplace} we conclude that
\begin{eqnarray*}
0 &\leq &\mathbf{E}\left[ \exp \left\{ -(1+o(1))W_{1,i-1,i}Q_{m}^{(i,N)}(%
\mathbf{s})\right\} \right] \\
&&-\mathbf{E}\left[ \exp \left\{ -(1+o(1))W_{1,i-1,i}Q_{m}^{(i,N)}(\mathbf{s}%
)-O(n^{-2\gamma _{i}}W_{1,i-1})\right\} \right] \\
&&\quad \leq 1-\mathbf{E}\left[ \exp \left\{ -O(n^{-\gamma
_{i+1}}W_{1,i-1})\right\} \right] =O\left( \left( n^{-\gamma _{i+1}}\right)
^{1/2^{i-1}}\right) =O\left( n^{-\gamma _{2}}\right) .
\end{eqnarray*}
As a result on account of (\ref{list}) we have for $m\sim yn^{\gamma
_{i}}$
\begin{eqnarray*}
Q_{m}^{(1,N)}(\mathbf{s}) &=&1-H_{m}^{(1,N)}\left( \mathbf{1}%
^{(i-1)},s_{i},...,s_{N}\right) +o(n^{-\gamma _{1}}) \\
&=&1-\mathbf{E}\left[ \exp \left\{ -(1+o(1))W_{1,i-1,i}Q_{m}^{(i,N)}(\mathbf{%
s})\right\} \right] +o(n^{-\gamma _{1}}) \\
&=&D_{i-1}\left( Q_{m}^{(i,N)}(\mathbf{s})\right) ^{1/2^{i-1}}+o(n^{-\gamma
_{1}}),
\end{eqnarray*}%
as required.

The lemma is proved.

\textbf{Proof of Lemma \ref{L_sharp!}}. Recalling Lemma \ref{L_initial} and
using the relation $\log (1-x)=-x+o(x),x\downarrow 0,$ we have
\begin{eqnarray}
&&\lim_{n\rightarrow \infty }\mathbf{E}^{n^{\gamma _{1}}}\left[ \exp \left\{
-\sum_{l=i}^{N}\lambda _{l}\frac{Z_{l}(m)}{n^{\left( l-i+1\right) \gamma
_{i}}}\right\} \right]   \notag \\
&&\quad =\lim_{n\rightarrow \infty }\exp \left\{ n^{\gamma _{1}}\log \mathbf{%
E}\left[ \exp \left\{ -\sum_{l=i}^{N}\lambda _{l}\frac{Z_{l}(m)}{n^{\left(
l-i+1\right) \gamma _{i}}}\right\} \right] \right\}   \notag \\
&&\quad =\exp \left\{ -\lim_{n\rightarrow \infty }n^{\gamma _{1}}\mathbf{E}%
\left[ 1-\exp \left\{ -\sum_{l=i}^{N}\lambda _{l}\frac{Z_{l}(m)}{n^{\left(
l-i+1\right) \gamma _{i}}}\right\} \right] \right\}   \notag \\
&&\quad =\exp \left\{ -D_{i-1}\left( \frac{\Phi _{i}(\lambda _{i}y,\lambda
_{i+1}y^{2},...,\lambda _{N}y^{N-i+1})}{y}\right) ^{1/2^{i-1}}\right\} .
\label{Analit1}
\end{eqnarray}%
Since the prelimiting function in $N-i+1$ complex variables $\lambda
_{i},...,\lambda _{N}$ is analytical and bounded in the domain $\left\{
Re\lambda _{l}>0,l=i,i+1,...,N\right\} :$
\begin{equation*}
\left\vert \mathbf{E}^{n^{\gamma _{1}}}\left[ \exp \left\{
-\sum_{l=i}^{N}\lambda _{l}\frac{Z_{l}(m)}{n^{\left( l-i+1\right) \gamma
_{i}}}\right\} \right] \right\vert \leq 1
\end{equation*}%
and converges for the real-valued $\lambda _{l}>0,l=i,i+1,...,N,$ it follows
by the Vitali and Weierstrass theorems that the limiting function is
analytical in the domain $\left\{ Re\lambda _{l}>0,l=i,i+1,...,N\right\} $
and, in addition, the derivative of the prelimiting function with resect to
any variable converges to the respective derivative of the limiting
function. Hence, on account of the equality
\begin{equation*}
\lim_{n\rightarrow \infty }\mathbf{E}\left[ \exp \left\{
-\sum_{l=i}^{N}\lambda _{l}\frac{Z_{l}(m)}{n^{\left( l-i+1\right) \gamma
_{i}}}\right\} \right] =1
\end{equation*}%
it is not difficult to deduce that
\begin{eqnarray*}
&&\lim_{n\rightarrow \infty }\frac{n^{\gamma _{1}}}{n^{\left( j-i+1\right)
\gamma _{i}}}\mathbf{E}\left[ Z_{j}(m)\exp \left\{ -\sum_{l=i}^{N}\lambda
_{l}\frac{Z_{l}(m)}{n^{\left( l-i+1\right) \gamma _{i}}}\right\} \right]
\mathbf{E}^{n^{\gamma _{1}}}\left[ \exp \left\{ -\sum_{l=i}^{N}\lambda _{l}%
\frac{Z_{l}(m)}{n^{\left( l-i+1\right) \gamma _{i}}}\right\} \right]  \\
&&\qquad \qquad \qquad \qquad \qquad =-\frac{\partial }{\partial \lambda _{j}%
}\exp \left\{ -D_{i-1}\left( \frac{\Phi _{i}(\lambda _{i}y,\lambda
_{i+1}y^{2},...,\lambda _{N}y^{N-i+1})}{y}\right) ^{1/2^{i-1}}\right\} ,
\end{eqnarray*}%
or, in view of (\ref{Analit1})
\begin{eqnarray*}
&&\lim_{n\rightarrow \infty }\frac{n^{\gamma _{1}}}{n^{\left( j-i+1\right)
\gamma _{i}}}\mathbf{E}\left[ Z_{j}(m)\exp \left\{ -\sum_{l=i}^{N}\lambda
_{l}\frac{Z_{l}(m)}{n^{\left( l-i+1\right) \gamma _{i}}}\right\} \right]  \\
&&\qquad \qquad \qquad \qquad \qquad \qquad \qquad =D_{i-1}\frac{\partial }{%
\partial \lambda _{j}}\left( \frac{\Phi _{i}(\lambda _{i}y,\lambda
_{i+1}y^{2},...,\lambda _{N}y^{N-i+1})}{y}\right) ^{1/2^{i-1}}.
\end{eqnarray*}%
The first part of Lemma \ref{L_sharp!} is proved.

To prove the second part it is necessary, basing on the representation
\begin{eqnarray*}
&&\mathbf{E}^{n^{\gamma _{1}}}\left[ 1-\left( 1-\exp \left\{
-\sum_{l=i}^{N}\lambda _{l}\frac{Z_{l}(m)}{n^{\left( l-i+1\right) \gamma
_{i}}}\right\} \right) I_{i-1}(m)\right]  \\
&=&(1+o(1))\exp \left\{ n^{\gamma _{1}}\mathbf{E}\left[ \left( 1-\exp
\left\{ -\sum_{l=i}^{N}\lambda _{l}\frac{Z_{l}(m)}{n^{\left( l-i+1\right)
\gamma _{i}}}\right\} \right) I_{i-1}(m)\right] \right\} ,
\end{eqnarray*}%
to repeat almost literally the arguments used earlier.

The lemma is proved.

\subsection{The case $n^{\protect\gamma _{i}}\ll m\ll n^{\protect\gamma %
_{i+1}},1\leq i\leq N-1$}

The aim of the present subsection is to check the validity of the following
statement:

\begin{lemma}
\label{L_interm!} If $n^{\gamma _{i}}\ll m\ll n^{\gamma _{i+1}}$ for some $%
i\in \left\{ 1,2,...,N-1\right\} ,$ then, for any $j\in \left\{
i,...,N-1\right\} $ and $\lambda _{l}\geq 0,l=i,...,N$
\begin{eqnarray*}
&&\lim_{n\rightarrow \infty }\frac{n^{\gamma _{1}}}{n^{\gamma
_{i+1}}m^{j-i-1}}\mathbf{E}\left[ Z_{j}(m)\exp \left\{
-\sum_{l=i+1}^{N}\lambda _{l}\frac{Z_{l}(m)}{n^{\gamma _{i+1}}m^{l-i-1}}%
\right\} I_{i}(m)\right] \\
&&\qquad =\lim_{n\rightarrow \infty }\frac{n^{\gamma _{1}}}{n^{\gamma
_{i+1}}m^{j-i-1}}\mathbf{E}\left[ Z_{j}(m)\exp \left\{
-\sum_{l=i+1}^{N}\lambda _{l}\frac{Z_{l}(m)}{n^{\gamma _{i+1}}m^{l-i-1}}%
\right\} \right] \\
&&\qquad =\frac{D_{i}a_{i+1,j}}{2^{i}}\left( \sum_{l=i+1}^{N}\lambda
_{l}a_{i+1,l}\right) ^{-1+1/2^{i}}.
\end{eqnarray*}
\end{lemma}

The needed result will be a corollary of a number of auxiliary statements.

\begin{lemma}
\label{L_intermONE}If $n^{\gamma _{i}}\ll m\ll n^{\gamma _{i+1}}$for some $%
i\in \left\{ 1,2,...,N-1\right\} $, then, for $j\geq i+1$ and $\lambda_l\geq 0,\, l=j,...,N$
\begin{equation*}
\lim_{n\rightarrow \infty }n^{\gamma _{i+1}}m^{j-i-1}\left( 1-\mathbf{E}_{j}%
\left[ \exp \left\{ -\sum_{l=j}^{N}\lambda _{l}\frac{Z_{l}(m)}{n^{\gamma
_{i+1}}m^{l-i-1}}\right\} \right] \right) =\sum_{l=j}^{N}\lambda _{l}a_{j,l}.
\end{equation*}
\end{lemma}

\textbf{Proof}. Put
\begin{equation}
s_{l}=\exp \left\{ -\frac{\lambda _{l}}{n^{\gamma _{i+1}}m^{l-i-1}}\right\}
,\ l=i+1,...,N.  \label{ss}
\end{equation}%
It is not difficult to check that, under the choice of variables
\begin{eqnarray*}
&&0\leq \sum_{l=j}^{N}\lambda _{l}\frac{\mathbf{E}_{j}\left[ Z_{l}(m)\right]
}{m^{l-j}}-n^{\gamma _{i+1}}m^{j-i-1}\left( 1-\mathbf{E}_{j}\left[
\prod_{l=j}^{N}s_{l}^{Z_{l}(m)}\right] \right) \\
&&\qquad \qquad \qquad \qquad \quad \leq \frac{m^{j-i-1}}{n^{\gamma _{i+1}}}%
\sum_{p,q=j}^{N}\lambda _{p}\lambda _{q}\frac{\mathbf{E}_{j}\left[
Z_{p}(m)Z_{q}(m)\right] }{m^{p-i-1}m^{q-i-1}}.
\end{eqnarray*}%
Using this inequality and the relations
\begin{equation*}
\sum_{l=j}^{N}\lambda _{l}\frac{\mathbf{E}_{j}\left[ Z_{l}(m)\right] }{%
m^{l-j}}\sim \sum_{l=j}^{N}\lambda _{l}a_{j,l}
\end{equation*}%
and
\begin{eqnarray*}
\sum_{p,q=j}^{N}\lambda _{p}\lambda _{q}\frac{\mathbf{E}_{j}\left[
Z_{p}(m)Z_{q}(m)\right] }{n^{2\gamma _{i+1}}m^{p-i-1}m^{q-i-1}} &\leq &\frac{%
C}{n^{2\gamma _{i+1}}}\sum_{p,q=j+2}^{N}\frac{m^{p+q-2j+1}}{%
m^{p-i-1}m^{q-i-1}} \\
&\leq &\frac{CN^{2}}{n^{2\gamma _{i+1}}}m^{2\left( i-j\right) +3}=\frac{%
CN^{2}}{n^{\gamma _{i+1}}m^{j-i-1}}\times \frac{m^{i-j+2}}{n^{\gamma _{i+1}}}
\\
&\leq &\frac{CN^{2}}{n^{\gamma _{i+1}}m^{j-i-1}}\times \frac{m}{n^{\gamma
_{i+1}}}=o\left( \frac{1}{n^{\gamma _{i+1}}m^{j-i-1}}\right) ,
\end{eqnarray*}%
following from Theorem \ref{T_Foster}, it is not difficult to demonstrate
the validity of the lemma.

\begin{lemma}
\label{L_interm2}If $n^{\gamma _{i}}\ll m\ll n^{\gamma _{i+1}}$ for some $%
i\in \left\{ 1,...,N-1\right\} $, then for $\lambda_l\geq 0,\, l=i+1,...,N$
\begin{eqnarray*}
&&\lim_{n\rightarrow \infty }n^{\gamma _{1}}\mathbf{E}\left[ \left( 1-\exp
\left\{ -\sum_{l=i+1}^{N}\lambda _{l}\frac{Z_{l}(m)}{n^{\gamma
_{i+1}}m^{l-i-1}}\right\} \right) I_{i}\left( m\right) \right] \\
&&\qquad =\lim_{n\rightarrow \infty }n^{\gamma _{1}}\mathbf{E}\left[ 1-\exp
\left\{ -\sum_{l=i+1}^{N}\lambda _{l}\frac{Z_{l}(m)}{n^{\gamma
_{i+1}}m^{l-i-1}}\right\} \right] \\
&&\qquad =D_{i}\left( \sum_{l=i+1}^{N}\lambda _{l}a_{i+1,l}\right)
^{1/2^{i}}.
\end{eqnarray*}
\end{lemma}

\textbf{Proof}. As before, it is sufficient to show the validity of the
second equality. Similarly to the arguments used earlier in the proof of
Lemma \ref{L_initial}, we have
\begin{equation*}
\mathbf{P}(T_{i}>m)\sim c_{1,i}m^{-2^{-(i-1)}}=o(n^{-\gamma _{1}}).
\end{equation*}%
Thus,
\begin{eqnarray*}
Q_{m}^{(1,N)}(\mathbf{s}) &=&\mathbf{E}\left[ 1-%
\prod_{j=1}^{N}s_{j}^{Z_{j}(m)}\right] \\
&=&\mathbf{E}\left[ \left( 1-\prod_{j=i+1}^{N}s_{j}^{Z_{j}(m)}\right)
;T_{i}\leq m\right] +o(n^{-\gamma _{1}}) \\
&=&1-H_{m}^{(1,N)}\left( \mathbf{1}^{(i)},s_{i+1},...,s_{N}\right)
+o(n^{-\gamma _{1}}).
\end{eqnarray*}%
Further,
\begin{eqnarray*}
&&H_{m}^{(1,N)}\left( \mathbf{1}^{(i)},s_{i+1},...,s_{N}\right) \\
&&\quad =\mathbf{E}\left[ \prod_{k=0}^{m-1}\prod_{r=1}^{i}%
\prod_{l=1}^{Z_{r}(k)}\prod_{j=i+1}^{N}\left( H_{m-k}^{(j,N)}(\mathbf{s}%
)\right) ^{\eta _{r,j}\left( k,l\right) }\right] \\
&&\quad =\mathbf{E}\left[ \prod_{k=0}^{m-1}\prod_{r=1}^{i}%
\prod_{l=1}^{Z_{r}(k)}\prod_{j=i+1}^{N}\left( H_{m-k}^{(j,N)}(\mathbf{s}%
)\right) ^{\eta _{r,j}\left( k,l\right) };T_{i}\leq \sqrt{mn^{\gamma _{i}}}%
\right] \\
&&\qquad +O\left( \mathbf{P}\left( T_{i}>\sqrt{mn^{\gamma _{i}}}\right)
\right) .
\end{eqnarray*}

Observing that $\lim_{m\rightarrow \infty }H_{m-k}^{(j,N)}(\mathbf{s})=1$
for $k\leq T_{i}\leq \sqrt{mn^{\gamma _{i}}}=o(m)$ and $j\geq i+1,$ we
conclude that, on the set $T_{i}\leq \sqrt{mn^{\gamma _{i}}}$%
\begin{eqnarray*}
&&\prod_{k=0}^{m-1}\prod_{r=1}^{i}\prod_{l=1}^{Z_{r}(k)}\prod_{j=i+1}^{N}%
\left( H_{m-k}^{(j,N)}(\mathbf{s})\right) ^{\eta _{r,j}\left( k,l\right) } \\
&&\quad =\exp \left\{
-\sum_{r=1}^{i}\sum_{k=0}^{T_{i}}\sum_{l=1}^{Z_{r}(k)}\sum_{j=i+1}^{N}\eta
_{r,j}\left( k,l\right) Q_{m-k}^{(j,N)}(\mathbf{s})(1+o(1))\right\} .
\end{eqnarray*}%
Lemma \ref{L_intermONE} and the estimate $m\ll n^{\gamma _{i+1}}$ give for $%
k=o\left( m\right) $ and $s_{l},l=i+1,...,N,$ from (\ref{ss}) :
\begin{equation}
Q_{m-k}^{(j,N)}(\mathbf{s})\sim Q_{m}^{(j,N)}(\mathbf{s})\sim \frac{1}{%
n^{\gamma _{i+1}}m^{j-i-1}}\sum_{l=j}^{N}\lambda _{l}a_{j,l}.  \label{tig}
\end{equation}%
Hence it follows that the relations
\begin{eqnarray*}
&&\sum_{r=1}^{i}\sum_{k=0}^{T_{i}}\sum_{l=1}^{Z_{r}(k)}\sum_{j=i+1}^{N}\eta
_{r,j}\left( k,l\right) Q_{m-k}^{(j,N)}(\mathbf{s}) \\
&&\quad =(1+o(1))\sum_{j=i+1}^{N}Q_{m}^{(j,N)}(\mathbf{s})\sum_{r=1}^{i}%
\sum_{k=0}^{T_{i}}\sum_{l=1}^{Z_{r}(k)}\eta _{r,j}\left( k,l\right) \\
&&\quad =(1+o(1))\sum_{j=i+1}^{N}W_{1,i,j}Q_{m}^{(j,N)}(\mathbf{s}) \\
&&\quad =(1+o(1))W_{1,i,i+1}Q_{m}^{(i+1,N)}(\mathbf{s})+O\left(
Q_{m}^{(i+2,N)}(\mathbf{s})W_{1,i}\right)
\end{eqnarray*}%
are valid on the set $T_{i}\leq \sqrt{mn^{\gamma _{i}}}=o(m)=o(n^{\gamma
_{i+1}}).$ Using the estimates
\begin{eqnarray*}
0 &\leq &\mathbf{E}\left[ \exp \left\{ -(1+o(1))W_{1,i,i+1}Q_{m}^{(i+1,N)}(%
\mathbf{s})\right\} \right] \\
&&-\mathbf{E}\left[ \exp \left\{ -(1+o(1))W_{1,i,i+1}Q_{m}^{(i+1,N)}(\mathbf{%
s})-O\left( Q_{m}^{(i+2,N)}(\mathbf{s})W_{1,i}\right) \right\} \right] \\
&&\qquad \qquad \leq 1-\mathbf{E}\left[ \exp \left\{ -O\left(
Q_{m}^{(i+2,N)}(\mathbf{s})W_{1,i}\right) \right\} \right] =O\left( \left(
\frac{1}{n^{\gamma _{i+1}}m}\right) ^{1/2^{i}}\right) \\
&&\qquad \qquad =o\left( \left( \frac{1}{n^{\gamma _{i+1}}n^{\gamma _{i}}}%
\right) ^{1/2^{i}}\right) =o\left( n^{-3\gamma _{1}/2}\right) =o\left(
n^{-\gamma _{1}}\right) ,
\end{eqnarray*}%
following from (\ref{Tot2}) and (\ref{tig}), and recalling (\ref{Tot1}) we
conclude that
\begin{eqnarray*}
Q_{m}^{(1,N)}(\mathbf{s}) &=&1-\exp \left\{ -(1+o(1))W_{1,i,i+1}n^{-\gamma
_{i+1}}\sum_{l=i+1}^{N}\lambda _{l}a_{i+1,l}\right\} +o\left( n^{-\gamma
_{1}}\right) \\
&=&D_{i}\left( \sum_{l=i+1}^{N}\lambda _{l}a_{i+1,l}\right)
^{1/2^{i}}n^{-\gamma _{1}}+o\left( n^{-\gamma _{1}}\right) ,
\end{eqnarray*}%
as required.

\textbf{Proof of Lemma \ref{L_interm!}}. To demonstrate the validity of
Lemma \ref{L_interm!} it is sufficient to recall Lemma \ref{L_interm2} and
to repeat (with evident changes) the arguments used to prove Lemma~\ref%
{L_sharp!}.

\subsection{The case $m\ll n^{\protect\gamma _{1}}$}

\begin{lemma}
\label{L_init!}If the parameters $m$ and $n$ tend to infinity in such a way
that $m\ll n^{\gamma _{1}},$ then for any $j\in \{1,...,N\}$ and $\lambda_l\geq 0,\, l=1,...,N$
\begin{equation*}
\lim_{m\rightarrow \infty }\frac{1}{m^{j-1}}\mathbf{E}\left[ Z_{j}(m)\exp
\left\{ -\sum_{l=1}^{N}\lambda _{l}\frac{Z_{l}(m)}{m^{l}}\right\} \right] =%
\frac{\partial \Phi _{1}(\lambda _{1},\lambda _{2},...,\lambda _{N})}{%
\partial \lambda _{j}}.
\end{equation*}
\end{lemma}

\textbf{Proof. }Recalling (\ref{Fost1}) and repeating the arguments used to
demonstrate Lemma \ref{L_sharp!}, we see that
\begin{eqnarray*}
Q_{m}^{(1,N)}(\mathbf{s}) &=&1-\exp \left\{ -(1+o(1))W_{1,i,i+1}n^{-\gamma
_{i+1}}\sum_{l=i+1}^{N}\lambda _{l}a_{i+1,l}\right\} +o\left( n^{-\gamma
_{1}}\right)  \\
&=&D_{i}\left( \sum_{l=i+1}^{N}\lambda _{l}a_{i+1,l}\right)
^{1/2^{i}}n^{-\gamma _{1}}+o\left( n^{-\gamma _{1}}\right) ,
\end{eqnarray*}%
as required.

Note, that%
\begin{equation}
\frac{\partial \Phi _{1}(\lambda _{1},\lambda _{2},...,\lambda _{N})}{%
\partial \lambda _{1}}\left\vert _{\mathbf{\lambda }=\mathbf{0}}\right. =1.
\label{NoAtom}
\end{equation}

\section{Proof of the limit theorems}

For $m<n$ introduce the functions
\begin{equation*}
\Psi ^{(i,N)}(m,n;\mathbf{s})=\mathbf{E}_{i}\left[ \mathbf{s}^{\mathbf{Z}%
(m)}I\left\{ T_{iN}=n\right\} \right] ,
\end{equation*}
where $I\left\{ A\right\} $ is the indicator of the event $A$. Our aim is to
investigate the asymptotic behavior of the quantity
\begin{equation*}
\mathbf{E}\left[ \mathbf{s}^{\mathbf{Z}(m)}|T_{1N}=n\right] =\frac{\Psi
^{(i,N)}(m,n;\mathbf{s})}{\mathbf{P}\left( T_{1N}=n\right) }
\end{equation*}
depending on the rate of growth of $n$ and $m$ to infinity. Clearly,
\begin{eqnarray*}
&&\Psi ^{(1,N)}(m,n;\mathbf{s}) \\
&=&\mathbf{E}\left[ \mathbf{s}^{\mathbf{Z}(m)}\left( \prod\limits_{l=1}^{N}%
\mathbf{P}_{l}^{Z_{l}(m)}\left( \mathbf{Z}(n-m)=\mathbf{0}\right)
-\prod\limits_{l=1}^{N}\mathbf{P}_{l}^{Z_{l}(m)}\left( \mathbf{Z}(n-m-1)=%
\mathbf{0}\right) \right) \right] .
\end{eqnarray*}
Using the formula
\begin{equation*}
\prod\limits_{l=1}^{N}X_{l}-\prod\limits_{l=1}^{N}Y_{l}=\sum_{j=1}^{N}\left(
X_{j}-Y_{j}\right)
\prod\limits_{l=1}^{j-1}Y_{l}\prod\limits_{l=j+1}^{N}X_{l},
\end{equation*}
where $\prod\limits_{l=1}^{0}Y_{l}=\prod\limits_{N+1}^{N}X_{l}=1,$ and
setting
\begin{equation*}
X_{l}=\mathbf{P}_{l}^{Z_{l}(m)}\left( \mathbf{Z}(n-m)=\mathbf{0}\right)
,\quad Y_{l}=\mathbf{P}_{l}^{Z_{l}(m)}\left( \mathbf{Z}(n-m-1)=\mathbf{0}%
\right) ,
\end{equation*}
we obtain
\begin{equation}
\Psi ^{(1,N)}(m,n;\mathbf{s})=\sum_{j=1}^{N}G_{j}\left( m,n;\mathbf{s}%
\right) ,  \label{dop1}
\end{equation}
where
\begin{equation}
G_{j}\left( m,n;\mathbf{s}\right) =\mathbf{E}\left[ \mathbf{s}^{\mathbf{Z}%
(m)}\left( X_{j}-Y_{j}\right)
\prod\limits_{l=1}^{j-1}Y_{l}\prod\limits_{l=j+1}^{N}X_{l}\right] .
\label{dop2}
\end{equation}
We separately investigate the behavior of the functions $G_{j}\left( m,n;%
\mathbf{s}\right) $ under an appropriate choice of the relationship between $%
m$ and $n$ and an appropriate choice of the components of $\mathbf{s}$.

Let%
\begin{equation}
x_{l}=\mathbf{P}_{l}\left( \mathbf{Z}(n-m)=\mathbf{0}\right) ,\quad y_{l}=%
\mathbf{P}_{l}\left( \mathbf{Z}(n-m-1)=\mathbf{0}\right) .  \label{SmalX}
\end{equation}
Then
\begin{eqnarray}
&&\mathbf{E}\left[ Z_{j}(m)\mathbf{s}^{\mathbf{Z}(m)}\left(
x_{j}-y_{j}\right)
y_{j}^{Z_{j}(m)}\prod\limits_{l=1}^{j-1}Y_{l}\prod\limits_{l=j+1}^{N}X_{l}%
\right]   \notag \\
&&\qquad \leq G_{j}\left( m,n;\mathbf{s}\right)   \notag \\
&&\qquad \quad \leq \mathbf{E}\left[ Z_{j}(m)\mathbf{s}^{\mathbf{Z}%
(m)}\left( x_{j}-y_{j}\right)
x_{j}^{Z_{j}(m)}\prod\limits_{l=1}^{j-1}Y_{l}\prod\limits_{l=j+1}^{N}X_{l}%
\right] .  \label{dop3}
\end{eqnarray}

In view of the asymptotic relations (\ref{SurvivSingle}) and Theorem \ref%
{T_loc}, we have
\begin{equation*}
\mathbf{P}_{i}\left( \mathbf{Z}(n)\neq \mathbf{0}\right) \sim \frac{c_{i,N}}{%
n^{\gamma _{i}}},\quad \mathbf{P}\left( T_{iN}=n\right) \sim \frac{g_{i,N}}{%
n^{1+\gamma _{i}}}.
\end{equation*}

Since
\begin{eqnarray*}
X_{l} &=&\exp \left\{ -Z_{l}(m)\mathbf{P}_{l}\left( \mathbf{Z}\left(
n-m\right) \neq \mathbf{0}\right) (1+\varepsilon _{l}\left( n,m\right)
)\right\} \\
&=&(1+\tilde{\varepsilon}_{l}\left( n,m\right) )\exp \left\{ -c_{l,N}\frac{%
Z_{l}(m)}{\left( n-m\right) ^{\gamma _{l}}}\right\} ,
\end{eqnarray*}
\begin{eqnarray*}
Y_{l} &=&\exp \left\{ -Z_{l}(m)\mathbf{P}_{l}\left( \mathbf{Z}\left(
n-m-1\right) \neq \mathbf{0}\right) (1+\varepsilon _{l}\left( n,m+1\right)
)\right\} \\
&=&(1+\tilde{\varepsilon}_{l}\left( n,m+1\right) )\exp \left\{ -c_{l,N}\frac{%
Z_{l}(m)}{\left( n-m\right) ^{\gamma _{l}}}\right\} ,
\end{eqnarray*}
it follows that for $m\ll n$ and $s_{l}=\exp \left\{ -\lambda
_{l}/L_{l}(m)\right\} $, where the functions $L_{l}(m),l=1,2,...,N$ will be
selected later on depending on the range of $m$ under consideration, it is
necessary to investigate, for each $j=1,2,...,N$ and up to negligible terms,
the asymptotic behavior of the quantity
\begin{eqnarray}
C_{j}(m,n) &=&\mathbf{E}\left[ Z_{j}(m)\mathbf{s}^{\mathbf{Z}(m)}\left(
x_{j}-y_{j}\right)
y_{j}^{Z_{j}(m)}\prod\limits_{l=1}^{j-1}X_{l}\prod\limits_{l=j+1}^{N}Y_{l}%
\right]  \notag \\
&=&\left( 1+\varepsilon _{j}(n,m)\right) \frac{g_{j,N}}{n^{1+\gamma _{j}}}%
\mathbf{E}\left[ Z_{j}(m)\mathbf{s}^{\mathbf{Z}(m)}y_{j}^{Z_{j}(m)}\prod%
\limits_{l=1}^{j-1}X_{l}\prod\limits_{l=j+1}^{N}Y_{l}\right]  \label{Evrist1}
\\
&=&\left( 1+\tilde{\varepsilon}_{j}(n,m)\right) \frac{g_{j,N}}{n^{1+\gamma
_{j}}}\mathbf{E}\left[ Z_{j}(m)\exp \left\{ -\sum_{l=1}^{N}\left( \frac{%
\lambda _{l}}{L_{l}(m)}+\frac{c_{l,N}}{\left( n-m\right) ^{\gamma _{l}}}%
\right) Z_{l}(m)\right\} \right] .  \notag
\end{eqnarray}

Consider first the case $m\ll n^{\gamma _{1}}$ and let $L_{l}(m)=m^{l}$.
Such a choice of parameters reduces (\ref{Evrist1}) to
\begin{equation}
C_{j}(m,n)=\left( 1+\varepsilon _{j}(n,m)\right) \frac{g_{j,N}}{n^{1+\gamma
_{j}}}\mathbf{E}\left[ Z_{j}(m)\exp \left\{ -\sum_{l=1}^{N}\lambda _{l}\frac{%
Z_{l}(m)}{m^{l}}\right\} \right] .  \label{Evrist2}
\end{equation}
These considerations lead to the following statement.

\begin{lemma}
\label{L_initial0}If $n^{\gamma _{1}}\gg m\rightarrow \infty $, then, for
all $\lambda _{l}\geq 0,l=1,...,N$
\begin{equation*}
\lim_{m\rightarrow \infty }\mathbf{E}\left[ \exp \left\{
-\sum_{l=1}^{N}\lambda _{l}\frac{Z_{l}(m)}{m^{l}}\right\} \Big|\,T_{N}=n%
\right] =\frac{\partial \Phi _{1}(\lambda _{1},\lambda _{2},...,\lambda _{N})%
}{\partial \lambda _{1}}.
\end{equation*}
\end{lemma}

\textbf{Proof}. We need to show that for $s_{l}=\exp \left\{ -\lambda
_{l}/m^{l}\right\} ,l=1,...,N,$
\begin{equation*}
\lim_{m\rightarrow \infty }\frac{\Psi ^{(1,N)}(m,n;\mathbf{s})}{\mathbf{P}%
\left( T_{N}=n\right) }=\frac{\partial \Phi _{1}(\lambda _{1},\lambda
_{2},...,\lambda _{N})}{\partial \lambda _{1}}.
\end{equation*}%
It follows from (\ref{dop1})--(\ref{Evrist2}), Theorem \ref{T_loc} and the
condition $m\ll n^{\gamma _{1}}$ that for every $j=1,2,...,N$ it is
necessary to investigate the asymptotic behavior of the quantity
\begin{equation*}
\frac{g_{j,N}}{g_{1,N}}\frac{n^{\gamma _{1}}}{n^{\gamma _{j}}}\mathbf{E}%
\left[ Z_{j}(m)\exp \left\{ -\sum_{l=1}^{N}\lambda _{l}\frac{Z_{l}(m)}{m^{l}}%
\right\} \right] .
\end{equation*}%
According to Lemma \ref{L_init!},
\begin{eqnarray*}
&&\lim_{m\rightarrow \infty }\frac{g_{j,N}}{g_{1,N}}\frac{n^{\gamma _{1}}}{%
n^{\gamma _{j}}}\mathbf{E}\left[ Z_{j}(m)\exp \left\{ -\sum_{l=1}^{N}\lambda
_{l}\frac{Z_{l}(m)}{m^{l}}\right\} \right] \\
&&\qquad \qquad \qquad =\frac{\partial \Phi _{1}(\lambda _{1},\lambda
_{2},...,\lambda _{N})}{\partial \lambda _{j}}\frac{g_{j,N}}{g_{1,N}}%
\lim_{m\rightarrow \infty }n^{\left( 1-2^{j-1}\right) \gamma _{1}}m^{j-1}.
\end{eqnarray*}%
Since
\begin{equation*}
\lim_{m\rightarrow \infty }n^{\left( 1-2^{j-1}\right) \gamma
_{1}}m^{j-1}=\delta _{1j},
\end{equation*}%
it follows that
\begin{equation*}
\lim_{m\rightarrow \infty }\frac{\Psi ^{(1,N)}(m,n;\mathbf{s})}{\mathbf{P}%
\left( T_{N}=n\right) }=\sum_{j=1}^{N}\lim_{m\rightarrow \infty }\frac{%
G_{j}\left( m,n;\mathbf{s}\right) }{\mathbf{P}\left( T_{N}=n\right) }=\frac{%
\partial \Phi _{1}(\lambda _{1},\lambda _{2},...,\lambda _{N})}{\partial
\lambda _{1}}.
\end{equation*}

The lemma is proved.

\begin{lemma}
\label{L_sharp00}If $m\sim yn^{\gamma _{i}}, y>0,$ for some $i\in \left\{
1,2,...,N-1\right\} ,$ then, for all $\lambda _{l}\geq 0,l=i,...,N$
\begin{eqnarray*}
&&\lim_{m\rightarrow \infty }\mathbf{E}\left[ \exp \left\{
-\sum_{l=i}^{N}\lambda _{l}\frac{Z_{l}(m)}{n^{\left( l-i+1\right) \gamma
_{i}}}\right\} I_{i-1}(m)\Big|\,T_{N}=n\right] \\
&&\quad =D_{i-1}\frac{g_{i,N}}{g_{1,N}}\frac{\partial }{\partial \lambda _{i}%
}\left( \frac{\Phi _{i}(\lambda _{i}^{\prime }y,\lambda _{i+1}^{\prime
}y^{2},\lambda _{i+2}y^{3},...,\lambda _{N}y^{N-i+1})}{y}\right) ^{1/2^{i-1}}
\\
&&\qquad +D_{i-1}\frac{g_{i+1,N}}{g_{1,N}}\frac{\partial }{\partial \lambda
_{i+1}}\left( \frac{\Phi _{i}(\lambda _{i}^{\prime }y,\lambda _{i+1}^{\prime
}y^{2},\lambda _{i+2}y^{3},...,\lambda _{N}y^{N-i+1})}{y}\right)
^{1/2^{i-1}},
\end{eqnarray*}%
where $\lambda _{i}^{\prime }=\lambda _{i}+c_{i,N},\lambda _{i+1}^{\prime
}=\lambda _{i+1}+c_{i+1,N}.$
\end{lemma}

\textbf{Proof}. Similarly to the proof of the previous lemma, it is
necessary to calculate, for each $j\in \left\{ i,i+1,...,N\right\} $ the
limit, as $m\rightarrow \infty $ of the quantity
\begin{eqnarray*}
&&\frac{g_{j,N}}{g_{1,N}}\frac{n^{\gamma _{1}}}{n^{\gamma _{j}}}\mathbf{E}%
\left[ Z_{j}(m)\exp \left\{ -\sum_{l=i}^{i+1}c_{l,N}\frac{Z_{l}(m)}{%
n^{(l-i+1)\gamma _{l}}}-\sum_{l=i}^{N}\lambda _{l}\frac{Z_{l}(m)}{n^{\left(
l-i+1\right) \gamma _{i}}}\right\} I_{i-1}(m)\right]  \\
&=&\frac{g_{j,N}}{g_{1,N}}\frac{n^{\left( j-i+1\right) \gamma _{i}}}{%
n^{\gamma _{j}}}\frac{n^{\gamma _{1}}}{n^{\left( j-i+1\right) \gamma _{i}}}%
\mathbf{E}\left[ Z_{j}(m)\exp \left\{ -\sum_{l=i}^{i+1}c_{l,N}\frac{Z_{l}(m)%
}{n^{(l-i+1)\gamma _{l}}}-\sum_{l=i}^{N}\lambda _{l}\frac{Z_{l}(m)}{%
n^{\left( l-i+1\right) \gamma _{i}}}\right\} I_{i-1}(m)\right]  \\
&=&(1+\varepsilon _{j}(n,m))D_{i-1}\frac{g_{j,N}}{g_{1,N}}\frac{n^{\left(
j-i+1\right) \gamma _{i}}}{n^{\gamma _{j}}}\frac{\partial }{\partial \lambda
_{j}}\left( \frac{\Phi _{i}(\lambda _{i}^{\prime }y,\lambda _{i+1}^{\prime
}y^{2},\lambda _{i+2}y^{3},...,\lambda _{N}y^{N-i+1})}{y}\right)
^{1/2^{i-1}},
\end{eqnarray*}%
where we have used Lemma \ref{L_sharp!} at the last step.

Hence the statement of the lemma follows easily, since
\begin{equation*}
\lim_{n\rightarrow \infty }\frac{n^{\left( j-i+1\right) \gamma _{i}}}{%
n^{\gamma _{j}}}=\left\{
\begin{array}{ccc}
1, & \text{if} & j=i,i+1, \\
&  &  \\
0, & \text{if} & j\neq i,i+1.%
\end{array}%
\right.
\end{equation*}

\begin{corollary}
\label{C_negl} Under conditions of Lemma \ref{L_sharp00}
\begin{equation}
\lim_{m\rightarrow \infty }\mathbf{P}(Z_{1}(m)+\cdots
+Z_{i-1}(m)>0|T_{N}=n)=0.  \label{kor1}
\end{equation}
\end{corollary}

\textbf{Proof.} Clearly,
\begin{eqnarray}
&&\mathbf{P}\left( Z_{1}(m)+\cdots +Z_{i-1}(m)=0|T_{N}=n\right) =\mathbf{E}%
\left[ I_{i-1}(m)|T_{N}=n\right]   \notag \\
&&\qquad\qquad\qquad\geq\mathbf{E}\left[ \exp \left\{ -\sum_{l=i}^{N}\lambda _{l}\frac{%
Z_{l}(m)}{n^{\left( l-i+1\right) \gamma _{i}}}\right\} I_{i-1}(m)\Big|%
\,T_{N}=n\right] .  \label{VV}
\end{eqnarray}%
Let $\mathbf{0}^{(N-i-1)}$ be an $N-i-1$--dimensional vector all whose
components are zeros. It follows from the definition of $\Phi _{i}$ (see (\ref{GlobalDif})) that%
\begin{equation*}
\frac{\Phi _{i}(\lambda _{i}y,\lambda _{i+1}y^{2},\mathbf{0}^{(N-i-1)})}{y}=%
\frac{\Phi _{i}^{\ast }(\lambda _{i}y,\lambda _{i+1}y^{2})}{y},
\end{equation*}%
where (recall (\ref{DefSimpl}))
\begin{equation*}
\frac{\Phi _{i}^{\ast }(\lambda _{i}y,\lambda _{i+1}y^{2})}{y}=\sqrt{\frac{%
m_{i,i+1}\lambda _{i+1}}{b_{i}}}\frac{b_{i}\lambda _{i}+\sqrt{%
b_{i}m_{i,i+1}\lambda _{i+1}}\tanh \left( y\sqrt{b_{i}m_{i,i+1}\lambda _{i+1}%
}\right) }{b_{i}\lambda _{i}\tanh \left( y\sqrt{b_{i}m_{i,i+1}\lambda _{i+1}}%
\right) +\sqrt{b_{i}m_{i,i+1}\lambda _{i+1}}}.
\end{equation*}%
Rather cumbersome calculations (which we omit), basing on the equalities,%
\begin{equation*}
c_{i,N}=\sqrt{b_{i}^{-1}m_{i,i+1}c_{i+1,N}},\;c_{1,N}=D_{i-1}\left(
c_{i,N}\right) ^{1/2^{i-1}},
\end{equation*}%
show that at the point $\left( \lambda _{i},\lambda _{i+1}\right) =\left(
c_{i,N},c_{i+1,N}\right) $
\begin{equation*}
D_{i-1}\frac{g_{i,N}}{g_{1,N}}\frac{\partial }{\partial \lambda _{i}}\left(
\frac{\Phi _{i}^{\ast }(\lambda _{i}y,\lambda _{i+1}y^{2})}{y}\right)
^{1/2^{i-1}}+D_{i-1}\frac{g_{i+1,N}}{g_{1,N}}\frac{\partial }{\partial
\lambda _{i+1}}\left( \frac{\Phi _{i}^{\ast }(\lambda _{i}y,\lambda
_{i+1}y^{2})}{y}\right) ^{1/2^{i-1}}=1.
\end{equation*}%
Combining this result with (\ref{VV}) and Lemma \ref{L_sharp00} gives%
\begin{equation*}
\lim \inf_{m\rightarrow \infty }\mathbf{P}\left( Z_{1}(m)+\cdots
+Z_{i-1}(m)=0|T_{N}=n\right) \geq 1.
\end{equation*}%
Thus,
\begin{equation*}
\lim_{m\rightarrow \infty }\mathbf{P}\left( Z_{1}(m)+\cdots
+Z_{i-1}(m)>0|T_{N}=n\right) =0.
\end{equation*}%
Corollary is proved.

\begin{lemma}
\label{L_intermed22}If $n^{\gamma _{i}}\ll m\ll n^{\gamma _{i+1}}$ for some $%
i\in \left\{ 1,2,...,N-1\right\} ,$ then, for any $\lambda _{l}\geq
0,l=i+1,...,N$%
\begin{eqnarray*}
&&\lim_{m\rightarrow \infty }\mathbf{E}\left[ \exp \left\{
-\sum_{l=i+1}^{N}\lambda _{l}\frac{Z_{l}(m)}{n^{\gamma _{i+1}}m^{l-i-1}}%
\right\} I_{i}(m)\Big|\,T_{N}=n\right]  \\
&&\qquad \qquad =\frac{D_{i}}{2^{i}}\frac{g_{i+1,N}}{g_{1,N}}\left(
c_{i+1,N}+\sum_{l=i+1}^{N}\lambda _{l}a_{i+1,l}\right) ^{-1+1/2^{i}}.
\end{eqnarray*}
\end{lemma}

\textbf{Proof}. Recalling (\ref{Evrist1}) and Lemma \ref{L_interm!}, we see
that it is necessary to calculate for $j\geq i+1$ the limit
\begin{eqnarray*}
&&\lim_{m\rightarrow \infty }\frac{g_{j,N}}{g_{1,N}}\frac{n^{\gamma _{1}}}{%
n^{\gamma _{j}}}\mathbf{E}\left[ Z_{j}(m)\exp \left\{ -c_{i+1,N}\frac{%
Z_{i+1}(m)}{n^{\gamma _{i+1}}}-\sum_{l=i+1}^{N}\lambda _{l}\frac{Z_{l}(m)}{%
n^{\gamma _{i+1}}m^{l-i-1}}\right\} I_{i}(m)\right] \\
&&\qquad =\lim_{m\rightarrow \infty }\frac{g_{j,N}}{g_{1,N}}\frac{n^{\gamma
_{i+1}}m^{j-i-1}}{n^{\gamma _{j}}}\frac{n^{\gamma _{1}}}{n^{\gamma
_{i+1}}m^{j-i-1}} \\
&&\qquad \qquad =\times \mathbf{E}\left[ Z_{j}(m)\exp \left\{ -c_{i+1,N}%
\frac{Z_{i+1}(m)}{n^{\gamma _{i+1}}}-\sum_{l=i+1}^{N}\lambda _{l}\frac{%
Z_{l}(m)}{n^{\gamma _{i+1}}m^{l-i-1}}\right\} I_{i}(m)\right] \\
&&\qquad =\frac{D_{i}a_{i+1,j}}{2^{i}}\frac{g_{j,N}}{g_{1,N}}\left(
c_{i+1,N}+\sum_{l=i+1}^{N}\lambda _{l}a_{i+1,l}\right)
^{-1+1/2^{i}}\lim_{m\rightarrow \infty }\frac{n^{\gamma _{i+1}}}{n^{\gamma
_{j}}}m^{j-i-1}.
\end{eqnarray*}%
If $j=i+1$, then
\begin{equation*}
\lim_{m\rightarrow \infty }\frac{n^{\gamma _{i+1}}}{n^{\gamma _{j}}}%
m^{j-i-1}=1.
\end{equation*}%
If $j>i+1$, then
\begin{equation*}
\lim_{m\rightarrow \infty }\frac{n^{\gamma _{i+1}}}{n^{\gamma _{j}}}%
m^{j-i-1}=0
\end{equation*}%
in view of the estimates
\begin{equation*}
\frac{n^{\gamma _{i+1}}}{n^{\gamma _{j}}}m^{j-i-1}\ll \frac{n^{\gamma
_{i+1}(j-i)}}{n^{\gamma _{j}}}=n^{\gamma _{i+1}\left( j-i-2^{j-i-1}\right)
}\leq 1.
\end{equation*}

Thus,
\begin{eqnarray*}
&&\lim_{m\rightarrow \infty }\mathbf{E}\left[ \exp \left\{
-\sum_{l=i+1}^{N}\lambda _{l}\frac{Z_{l}(m)}{n^{\gamma _{i+1}}m^{l-i-1}}%
\right\} I_{i}(m)\Big|\,T_{N}=n\right]  \\
&&\qquad \qquad \qquad =\frac{D_{i}}{2^{i}}\frac{g_{i+1,N}}{g_{1,N}}\left(
c_{i+1,N}+\sum_{l=i+1}^{N}\lambda _{l}a_{i+1,l}\right) ^{-1+1/2^{i}}.
\end{eqnarray*}

Lemma \ref{L_intermed22} is proved.

\begin{corollary}
\label{C_negl1}Under conditions of Lemma \ref{L_intermed22}
\begin{equation*}
\lim_{m\rightarrow \infty }\mathbf{P}(Z_{1}(m)+\cdots +Z_{i}(m)>0|T_{N}=n)=0.
\end{equation*}
\end{corollary}

\textbf{Proof.} In virtue of Lemma 6 in \cite{VV14} and the equalities $%
g_{k,N}=\gamma _{k}c_{k,N},\gamma _{i+1}=2^{i}\gamma _{1},$ we have
\begin{equation*}
\frac{D_{i}}{2^{i}}\left( c_{i+1,N}\right) ^{-1+1/2^{i}}\frac{g_{i+1,N}}{%
g_{1,N}}=\frac{D_{i}\left( c_{i+1,N}\right) ^{1/2^{i}}\gamma _{i+1}}{%
2^{i}\gamma _{1}c_{1,N}}=\frac{D_{i}\left( c_{i+1,N}\right) ^{1/2^{i}}}{%
c_{1,N}}=1,
\end{equation*}%
which in view of Lemma \ref{L_intermed22} finishes the proof.

\begin{lemma}
\label{L_final}If $m\sim xn,x\in \left( 0,1\right),$  then, for $\lambda _{N}\geq 0$
\begin{eqnarray*}
&&\lim_{n\rightarrow \infty }n^{\gamma _{1}}\mathbf{E}\left[ \left( 1-\exp
\left\{ -\lambda _{N}\frac{Z_{N}(m)}{b_{N}n}\right\} \right) I_{N-1}(m)%
\right]  \\
&&\qquad \qquad =\lim_{n\rightarrow \infty }n^{\gamma _{1}}\mathbf{E}\left[
1-\exp \left\{ -\lambda _{N}\frac{Z_{N}(m)}{b_{N}n}\right\} \right] =\frac{%
c_{1,N}}{x^{\gamma _{1}}}\left( 1-\frac{1}{1+x\lambda _{N}}\right) ^{\gamma
_{1}}.
\end{eqnarray*}
\end{lemma}

\textbf{Proof}. This statement follows from Theorem 4 in \cite{VV14} and the
asymptotic representation (\ref{SurvivSingle}).

\begin{lemma}
\label{L_final2}If $m\sim xn,x\in \left( 0,1\right),$ then for $\lambda _{N}\geq 0$
\begin{eqnarray*}
&&\lim_{m\rightarrow \infty }\mathbf{E}\left[ \exp \left\{ -\lambda _{N}%
\frac{Z_{N}(m)}{b_{N}n}\right\} I_{N-1}(m)\Big|\,T_{N}=n\right]  \\
&&\qquad \qquad =\lim_{m\rightarrow \infty }\mathbf{E}\left[ \exp \left\{
-\lambda _{N}\frac{Z_{N}(m)}{b_{N}n}\right\} \Big|\,T_{N}=n\right]  \\
&&\qquad \qquad =\frac{1}{\left( 1+\left( 1-x\right) \lambda _{N}\right)
^{1-\gamma _{1}}}\frac{1}{\left( 1+x\left( 1-x\right) \lambda _{N}\right)
^{1+\gamma _{1}}}.
\end{eqnarray*}
\end{lemma}

\textbf{Proof}. Similarly to the proof  of (\ref{Evrist1}) one can show,
using the notations from (\ref{SmalX}), that for $m\sim xn,x\in \left(
0,1\right) $
\begin{eqnarray*}
&&\mathbf{E}\left[ \exp \left\{ -\lambda _{N}\frac{Z_{N}(m)}{b_{N}n}\right\}
I_{N-1}(m);\,T_{N}=n\right]  \\
&=&\mathbf{E}\left[ \exp \left\{ -\lambda _{N}\frac{Z_{N}(m)}{b_{N}n}%
\right\} \left( x_{N}^{Z_{N}(m)}-y_{N}^{Z_{N}(m)}\right) I_{N-1}(m)\right]
\\
&=&\frac{(1+\varepsilon _{1}(m,n))g_{N,N}}{(n(1-x))^{2}}\mathbf{E}\left[
Z_{N}(m)\exp \left\{ -\left( \lambda _{N}+\frac{b_{N}c_{N,N}}{1-x}\right)
\frac{Z_{N}(m)}{b_{N}n}\right\} I_{N-1}(m)\right] .
\end{eqnarray*}%
By Lemma \ref{L_final} for $\lambda \geq 0$ we have%
\begin{eqnarray*}
&&\lim_{n\rightarrow \infty }\frac{n^{\gamma _{1}}}{b_{N}n}\mathbf{E}\left[
Z_{N}(m)\exp \left\{ -\lambda \frac{Z_{N}(m)}{b_{N}n}\right\} I_{N-1}(m)%
\right]  \\
&&\qquad =\lim_{n\rightarrow \infty }\frac{\partial }{\partial \lambda }%
n^{\gamma _{1}}\mathbf{E}\left[ \left( 1-\exp \left\{ -\lambda \frac{Z_{N}(m)%
}{b_{N}n}\right\} \right) I_{N-1}(m)\right]  \\
&&\qquad =\frac{\partial }{\partial \lambda }\lim_{n\rightarrow \infty
}n^{\gamma _{1}}\mathbf{E}\left[ \left( 1-\exp \left\{ -\lambda \frac{%
Z_{N}(m)}{b_{N}n}\right\} \right) I_{N-1}(m)\right]  \\
&&\qquad =\frac{\partial }{\partial \lambda }\frac{c_{1,N}}{x^{\gamma _{1}}}%
\left( 1-\frac{1}{1+x\lambda }\right) ^{\gamma _{1}}=\gamma
_{1}c_{1,N}\left( \frac{\lambda }{1+x\lambda }\right) ^{\gamma _{1}-1}\frac{1%
}{\left( 1+x\lambda \right) ^{2}}.
\end{eqnarray*}%
Hence, taking into account the equalities $g_{N,N}=b_{N}^{-1},$ $%
g_{1,N}=\gamma _{1}c_{1,N},$ $b_{N}c_{1,N}=1,$using the relation%
\begin{equation*}
\mathbf{P}\left( T_{N}=n\right) \sim \frac{g_{1,N}}{n^{1+\gamma _{1}}},
\end{equation*}%
and setting%
\begin{equation*}
\lambda =\lambda _{N}+\frac{b_{N}c_{N,N}}{1-x}=\lambda _{N}+\frac{1}{1-x},
\end{equation*}%
after evident simplifications we obtain
\begin{eqnarray*}
&&\lim_{n\rightarrow \infty }\frac{\gamma _{1}c_{1,N}b_{N}g_{N,N}n^{\gamma
_{1}+1}}{g_{1,N}(n(1-x))^{2}b_{N}}\mathbf{E}\left[ Z_{N}(m)\exp \left\{
-\left( \lambda _{N}+\frac{1}{1-x}\right) \frac{Z_{N}(m)}{b_{N}n}\right\}
I_{N-1}(m)\right]  \\
&&\qquad =\frac{1}{\left( 1+\left( 1-x\right) \lambda _{N}\right) ^{1-\gamma
_{1}}}\frac{1}{\left( 1+x\left( 1-x\right) \lambda _{N}\right) ^{1+\gamma
_{1}}},
\end{eqnarray*}%
as required.

\begin{corollary}
\label{C_negl2}Under conditions of Lemma \ref{L_final2}
\begin{equation*}
\lim_{m\rightarrow \infty }\mathbf{P}(Z_{1}(m)+\cdots
+Z_{N-1}(m)>0|T_{N}=n)=0.
\end{equation*}
\end{corollary}

\textbf{Proofs of Theorems \ref{T_initial}--\ref{T_finalstage}}. The
statements of Theorems \ref{T_initial}--\ref{T_finalstage} follow in an
evident way from Lemmas \ref{L_initial0}, \ref{L_sharp00}, \ref{L_intermed22}%
, \ref{L_final2} and Corollaries \ref{C_negl}, \ref{C_negl1}, \ref{C_negl2}.

\end{document}